\documentclass[11pt, a4paper]{article}
\usepackage{amsmath}
\usepackage{amsfonts}
\usepackage{amssymb}
\usepackage{ulem}

\usepackage{amsthm}
\usepackage{url}
\usepackage{dcolumn}

\usepackage{cite}
\usepackage{fancyhdr}

\usepackage{graphicx, color}
\usepackage{hyperref}
\usepackage{xcolor}

\newtheorem*{rem*}{Remark}

\begin{document}

\title{Convex pentagons and concave octagons that can form \\rotationally 
symmetric tilings}
\author{ Teruhisa SUGIMOTO$^{ 1), 2)}$ }
\date{}
\maketitle

{\footnotesize

\begin{center}
$^{1)}$ The Interdisciplinary Institute of Science, Technology and Art

$^{2)}$ Japan Tessellation Design Association

E-mail: ismsugi@gmail.com
\end{center}

}

{\small
\begin{abstract}
\noindent
In this study, the properties of convex pentagons that can form 
rotationally symmetric edge-to-edge tilings are discussed. Because the 
rotationally symmetric tilings are formed by concave octagons that 
are generated by two convex pentagons connected through a line symmetry, 
they are considered to be equivalent to rotationally symmetric tilings with 
concave octagons. In addition, under certain circumstances, tiling-like 
patterns with a regular polygonal hole at the center can be formed using 
these convex pentagons.
\end{abstract}
}


%%%%%%%%%%%%%%%%%%%%%%%%%%%%%%%%%%%%%%%%%%%%%%%%%%%%%%%%%%%%%%%%%%%%%%
%%%%%%%%%%%%%%%%%%%%%%%%%%%%%%%%%%%%%%%%%%%%%%%%%%%%%%%%%%%%%%%%%%%%%%
\section{Introduction}
\label{section1}

\begin{figure}[htbp]
 \centering\includegraphics[width=15.cm,clip]{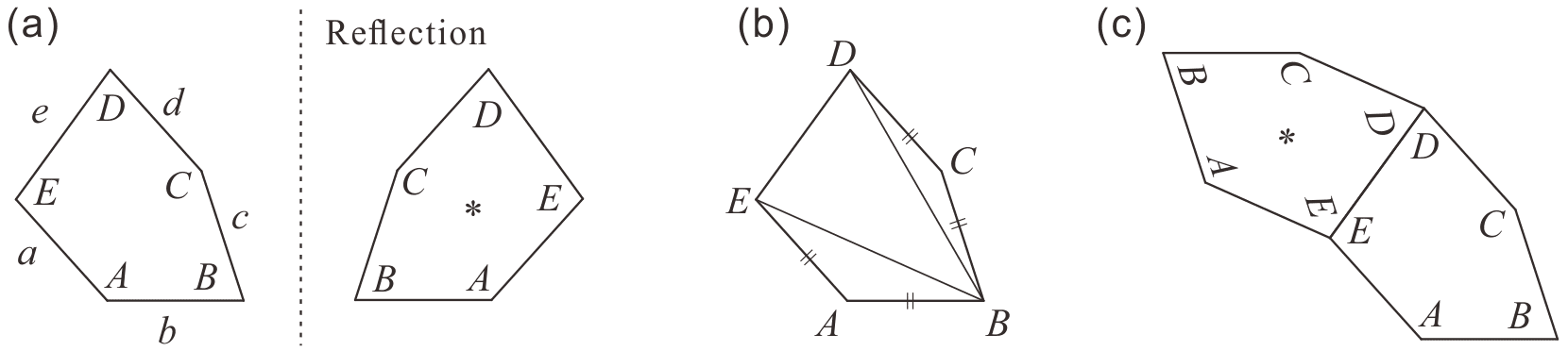} 
  \caption{{\small 
Nomenclature for vertices and edges of a convex pentagon,  
three triangles in the convex pentagonal tile C11-T1a, and the Octa-unit
} 
\label{fig01}
}
\end{figure}

In this study, as shown in Figure~\ref{fig01}(a), let us label the vertices (interior angles) 
of the convex pentagon $A$, $B$, $C$, $D$, and $E$, and its edges $a$, $b$, $c$, 
$d$, and $e$ in a fixed manner. A convex pentagonal tile\footnote{ 
A \textit{tiling} (or \textit{tessellation}) of the plane is a collection of sets that 
are called tiles, which covers a plane without gaps and overlaps, except for the 
boundaries of the tiles. The term ``tile" refers to a topological disk, whose boundary 
is a simple closed curve. If all the tiles in a tiling are of the same size and shape, 
then the tiling is \textit{monohedral}~\cite{G_and_S_1987, wiki_pentagon_tiling}. 
In this study, a polygon that admits a monohedral tiling is called a 
\textit{polygonal tile}~\cite{Sugimoto_NoteTP, Sugimoto_2012, Sugimoto_2016}. 
Note that, in monohedral tiling, it admits the use of reflected tiles.
} that satisfies the conditions ``$A + B + C = 360^ \circ ,\;C = 2D,\;a = b = c = d$'' 
belongs to the Type 1 family\footnote{ 
To date, fifteen families of convex pentagonal tiles, each of them referred to as 
a ``Type," are known~\cite{G_and_S_1987, Sugimoto_NoteTP, wiki_pentagon_tiling}. 
For example, if the sum of three consecutive angles in a convex pentagonal tile 
is  $360^ \circ $, the pentagonal tile belongs to the Type 1 family. Convex 
pentagonal tiles belonging to some families also exist. Known convex pentagonal 
tiles can form periodic tiling. In May 2017, Micha\"{e}l Rao declared that the 
complete list of Types of convex pentagonal tiles had been obtained (i.e., they 
have only the known 15 families), but it does not seem to be fixed as of March 
2020~\cite{wiki_pentagon_tiling}.}~\cite{G_and_S_1987, Sugimoto_2012, 
Sugimoto_NoteTP, wiki_pentagon_tiling
}.  The tile is called ``C11-T1A'' in \cite{S_and_O_2009} and its geometric properties 
are shown. Because this convex pentagon has four equal-length edges, its interior 
can be divided into a triangle \textit{BDE} and two isosceles triangles \textit{ABE} 
and \textit{BCD} as shown in Figure~\ref{fig01}(b). Given the foregoing properties, 
the relational expression of the interior angle of each vertex of C11-T1A can be 
rewritten as follows:

\renewcommand{\figurename}{{\small Figure.}}

\renewcommand{\figurename}{{\small Figure.}}
\begin{figure}[htbp]
 \centering\includegraphics[width=15cm,clip]{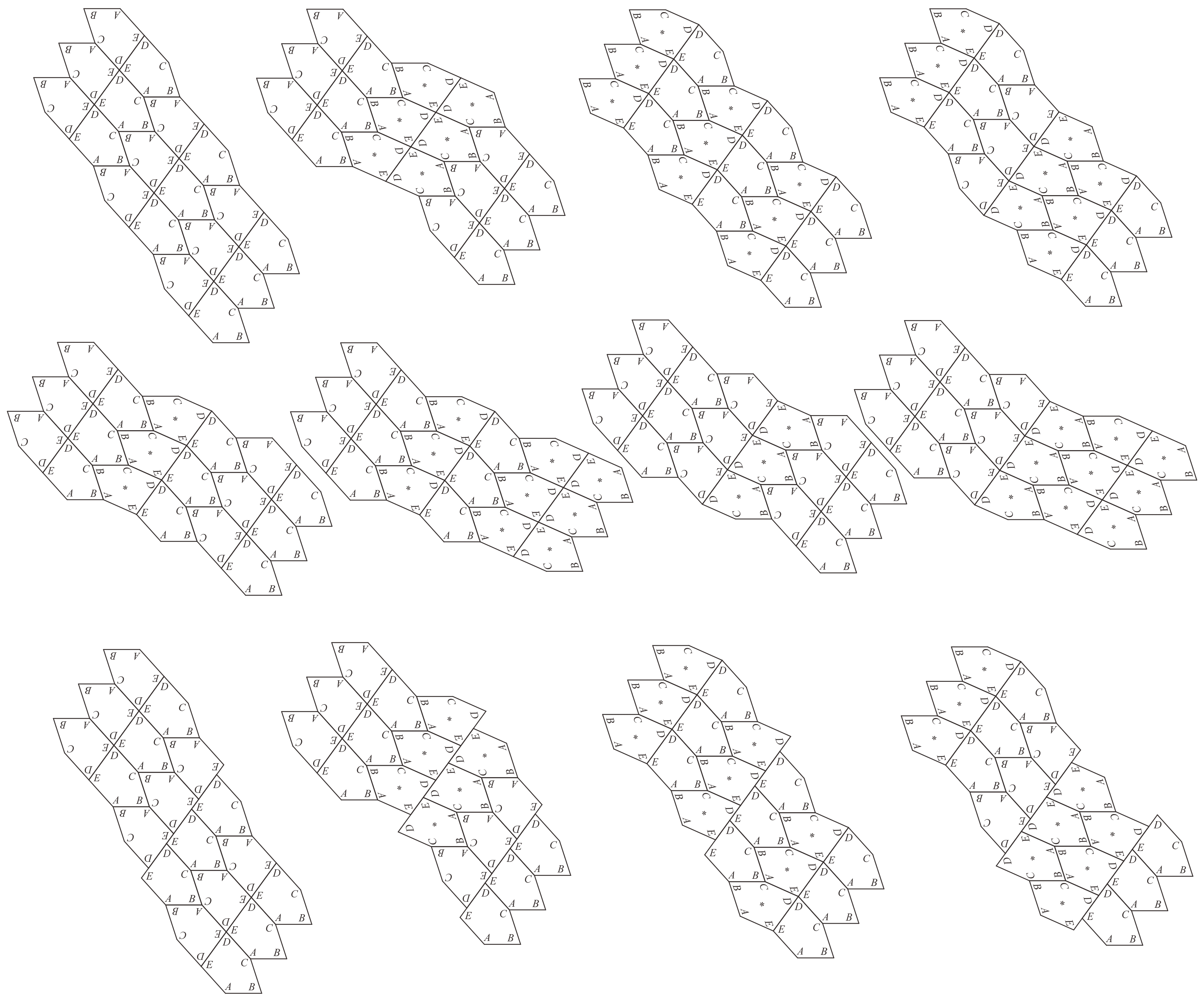} 
  \caption{{\small 
Examples of variations of Type 1 tilings by C11-T1A
} 
\label{fig02}
}
\end{figure}

\begin{equation}
\label{eq1}
\left\{ {\begin{array}{l}
 A = 180^ \circ - \dfrac{2B}{3\strut}, \\ 
 C = 180^ \circ - \dfrac{B}{3\strut}, \\ 
 D = 90^ \circ - \dfrac{B}{6\strut}, \\ 
 E = 90^ \circ + \dfrac{B}{6\strut}, \\ 
 \end{array}} \right.
\end{equation}

\noindent
where $0^ \circ < B < 180^ \circ $. From the relationship between the five 
interior angles, the vertex concentrations that are always valid in tilings are 
``$A + B + C = 360^ \circ ,\;2E + C = 360^ \circ ,\;2D + A + B = 360^ \circ ,\;
2D + 2E = 360^ \circ $." If $a = b = c = d = 1$, then the edge length of $e$ 
can be expressed as follows:

\[
e = 2\sin \left( {\frac{B}{2}} \right).
\]

\begin{figure}[htb]
 \centering\includegraphics[width=15cm,clip]{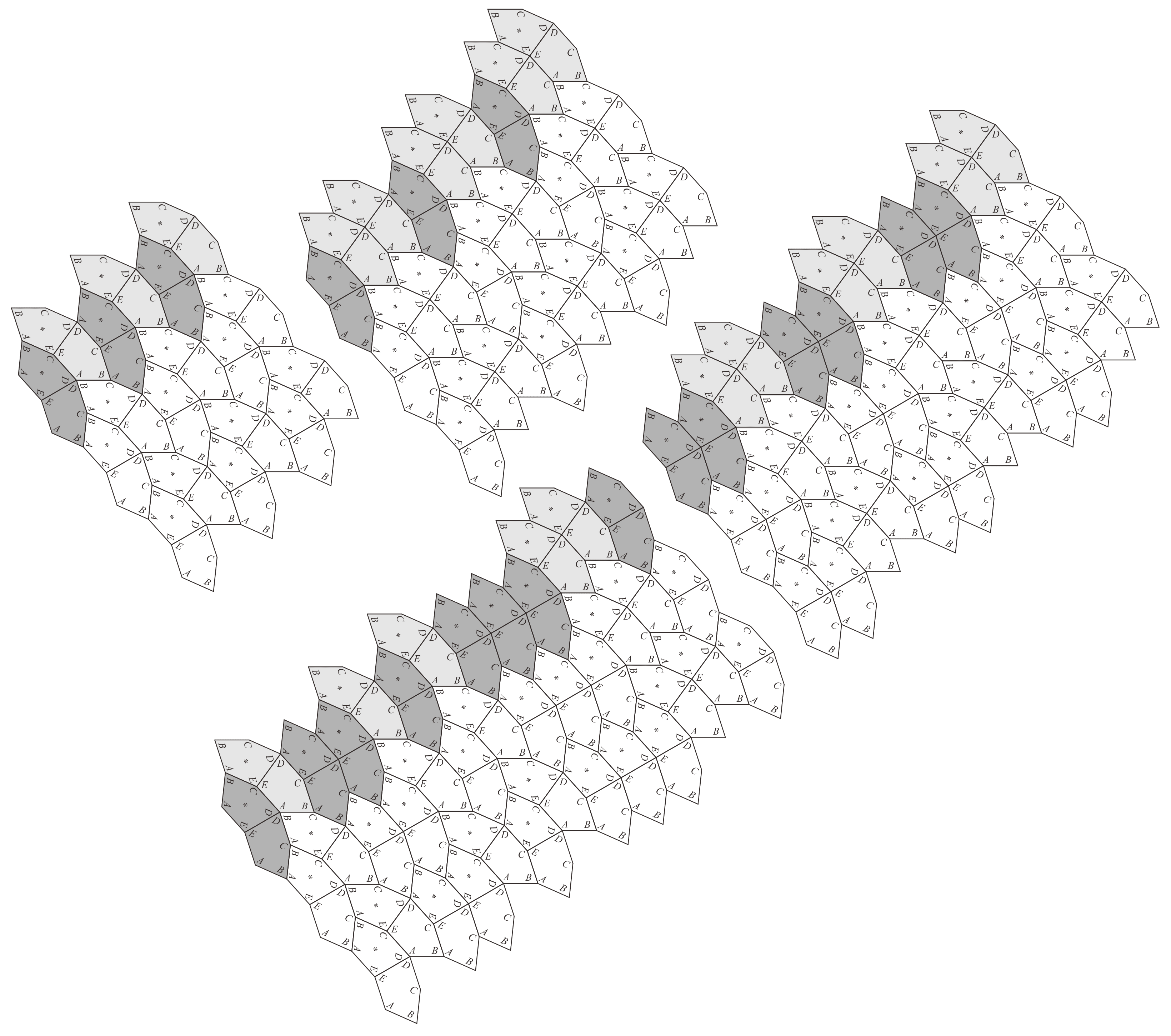} 
  \caption{{\small 
Examples of tilings with Octa-units of C11-T1A
} 
\label{fig03}
}
\end{figure}

As shown in Figure~\ref{fig02}, C11-T1A can form the representative tiling of Type 1 
or variations of Type 1 tilings (i.e., tilings whose vertices are formed only by the 
relations of $A + B + C = 360^ \circ $ and $D + E = 180^ \circ$). The edge $e$ of 
C11-T1A is the only edge of different length. As shown in Figure~\ref{fig01}(c), an 
equilateral concave octagon formed by two convex pentagons, connected through a 
line symmetry whose axis is edge $e$, is referred to as the \textit{Octa-unit}. 
C11-T1A can form tilings that are not contained in the variations of Type 1 
tilings, as shown in Figure~\ref{fig03}, by using Octa-units. That is, C11-T1A can 
also form tilings by freely combining Octa-units with different directions 
in one direction~\cite{Iliev_2018, S_and_O_2009}

In this study, we introduce the convex pentagon C11-T1A that can generate 
countless numbers of rotationally symmetric tilings. Because the Octa-unit of 
C11-T1A is used in these rotationally symmetric tilings, it is considered 
similar to the concave octagons that can generate countless rotationally 
symmetric tilings.

\renewcommand{\figurename}{{\small Figure.}}

%%%%%%%%%%%%%%%%%%%%%%%%%%%%%%%%%%%%%%%%%%%%%%%%%%%%%%%%%%%%%%%%%%%%%%
%%%%%%%%%%%%%%%%%%%%%%%%%%%%%%%%%%%%%%%%%%%%%%%%%%%%%%%%%%%%%%%%%%%%%%
\section{Rotationally symmetric tilings}
\label{section2}

In \cite{Iliev_2018}, there are figures of four-fold, five-fold, six-fold, and ten-fold 
rotationally symmetric edge-to-edge\footnote{ 
A tiling by convex polygons is 
\textit{edge-to-edge} if any two convex polygons in a tiling are either disjoint 
or share one vertex or an entire edge in common. Then other case is 
\textit{non-edge-to-edge}~\cite{G_and_S_1987, Sugimoto_2012, Sugimoto_NoteTP }.
} tilings using convex pentagonal tiles (or tiles that can be regarded as concave 
octagons). In \cite{Iliev_2018}, there are multiple types of edge-to-edge tilings with 
six-fold rotational symmetry by an equilateral convex pentagon, including 
Hirschhorn and Hunt's types~\cite{H_and_H_1985}. The convex pentagons forming 
these rotationally symmetric edge-to-edge tilings corresponded to C11-T1A.

By using the relational expression of the interior angles of (\ref{eq1}), the 
conditions of pentagonal tile C11-T1A that can form $n$-fold rotationally 
symmetric edge-to-edge tilings are expressed in (\ref{eq2}).

\begin{equation}
\label{eq2}
\left\{ {\begin{array}{l}
 A = 180^ \circ - \dfrac{240^ \circ}{n\strut}, \\ 
 B = \dfrac{360^ \circ}{n\strut}, \\ 
 C = 180^ \circ - \dfrac{120^ \circ}{n\strut}, \\ 
 D = 90^ \circ - \dfrac{60^ \circ}{n\strut}, \\ 
 E = 90^ \circ + \dfrac{60^ \circ}{n\strut}, \\ 
 a = b = c = d, \\ 
 \end{array}} \right.
\end{equation}

\noindent
where $n$ is an integer greater than or equal to three, because 
$0^ \circ < B < 180^ \circ $. 

Table~\ref{tab1} presents some of the relationships between the interior angles of 
convex pentagons satisfying (\ref{eq2}) that can form the $n$-fold rotationally 
symmetric edge-to-edge tilings. (For $n = 3\!-\!13$, tilings with convex pentagonal 
tiles are drawn. For further details, Figures~\ref{fig05}--\ref{fig15}.) Given 
that $n = 6$, the convex pentagonal tile that satisfies (\ref{eq2}) is an equilateral 
convex pentagon, and can form several different six-fold rotationally 
symmetric edge-to-edge tilings~\cite{H_and_H_1985, Iliev_2018, S_and_O_2009, 
wiki_pentagon_tiling}. The tiling in Figure~\ref{fig08} is one of the six-fold rotationally 
symmetric edge-to-edge tilings by the equilateral convex pentagon that 
Hirschhorn and Hunt presented~\cite{H_and_H_1985}. Note that the $n$-fold 
rotationally symmetric edge-to-edge tilings by convex pentagonal tiles 
satisfying (\ref{eq2}) have $C_{n}$ symmetry\footnote{ 
``$C_{n}$" is based on the Schoenflies notation for symmetry in a 
two-dimensional point group~\cite{wiki_point_group, wiki_schoenflies_notation}. 
``$C_{n}$" represents an $n$-fold rotation axis without reflection. The notation 
for symmetry is based on that presented in \cite{Klaassen_2016}.
} because they have rotational symmetry, but no axis of reflection symmetry.

\begin{table}[t]
 \begin{center}
{\small
\caption[Table 1]{
Interior angles of convex pentagons satisfying (\ref{eq2}) 
that can form the $n$-fold rotationally symmetric edge-to-edge tilings
}
\label{tab1}
}
\begin{tabular}
{c| D{.}{.}{2} D{.}{.}{2} D{.}{.}{2} D{.}{.}{2} D{.}{.}{2} |c|c}
%{c|rrrrr|c|c}
\hline
\raisebox{-1.50ex}[0cm][0cm]{$n$}& 
\multicolumn{5}{c|}{\shortstack{ Value of interior angle (degree) } } & 
\raisebox{-4.6ex}[0.7cm][0.5cm]{\shortstack{Edge \\length \\of $e$}} & 
\raisebox{-3.0ex}[0.7cm][0.5cm]{\shortstack{Figure \\number}} \\

 & 
\textit{A} & 
\textit{B}& 
\textit{C}& 
\textit{D}& 
\textit{E}& 
 & 
  \\
\hline
3& 
100& 
120& 
140& 
70& 
110& 
1.732 & 
\ref{fig05} \\
\hline
4& 
120& 
90& 
150& 
75& 
105& 
1.414 & 
\ref{fig06} \\
\hline
5& 
132& 
72& 
156& 
78& 
102& 
1.176 & 
\ref{fig07} \\
\hline
6& 
140& 
60& 
160& 
80& 
100& 
1 & 
\ref{fig08} \\
\hline
7& 
145.71 & 
51.43 & 
162.86 & 
81.43 & 
98.57 & 
0.868 & 
\ref{fig09} \\
\hline
8& 
150& 
45& 
165& 
82.5& 
97.5& 
0.765 & 
\ref{fig10} \\
\hline
9& 
153.33 & 
40 & 
166.67 & 
83.33 & 
96.67 & 
0.684 & 
\ref{fig11} \\
\hline
10& 
156& 
36& 
168& 
84& 
96& 
0.618 & 
\ref{fig12} \\
\hline
11& 
158.18 & 
32.73 & 
169.09 & 
84.55 & 
95.45 & 
0.563 & 
\ref{fig13} \\
\hline
12& 
160& 
30& 
170& 
85& 
95& 
0.518 & 
\ref{fig14} \\
\hline
13& 
161.54 & 
27.69 & 
170.77 & 
85.38 & 
94.62 & 
0.479 & 
\ref{fig15} \\
\hline
14& 
162.86 & 
25.71 & 
171.43 & 
85.71 & 
94.29 & 
0.445 & 
 \\
\hline
15& 
164& 
24& 
172& 
86& 
94& 
0.416 & 
 \\
\hline
16& 
165& 
22.5& 
172.5& 
86.25& 
93.75& 
0.390 & 
 \\
\hline
17& 
165.88 & 
21.18 & 
172.94 & 
86.47 & 
93.53 & 
0.367 & 
 \\
\hline
18& 
166.67 & 
20 & 
173.33 & 
86.67 & 
93.33 & 
0.347 &
 \\
\hline
...& 
...& 
...& 
...& 
...& 
...& 
...& 
 \\
\hline
\end{tabular}

\end{center}

\end{table}

\renewcommand{\figurename}{{\small Figure.}}
\begin{figure}[!h]
 \centering\includegraphics[width=14.5cm,clip]{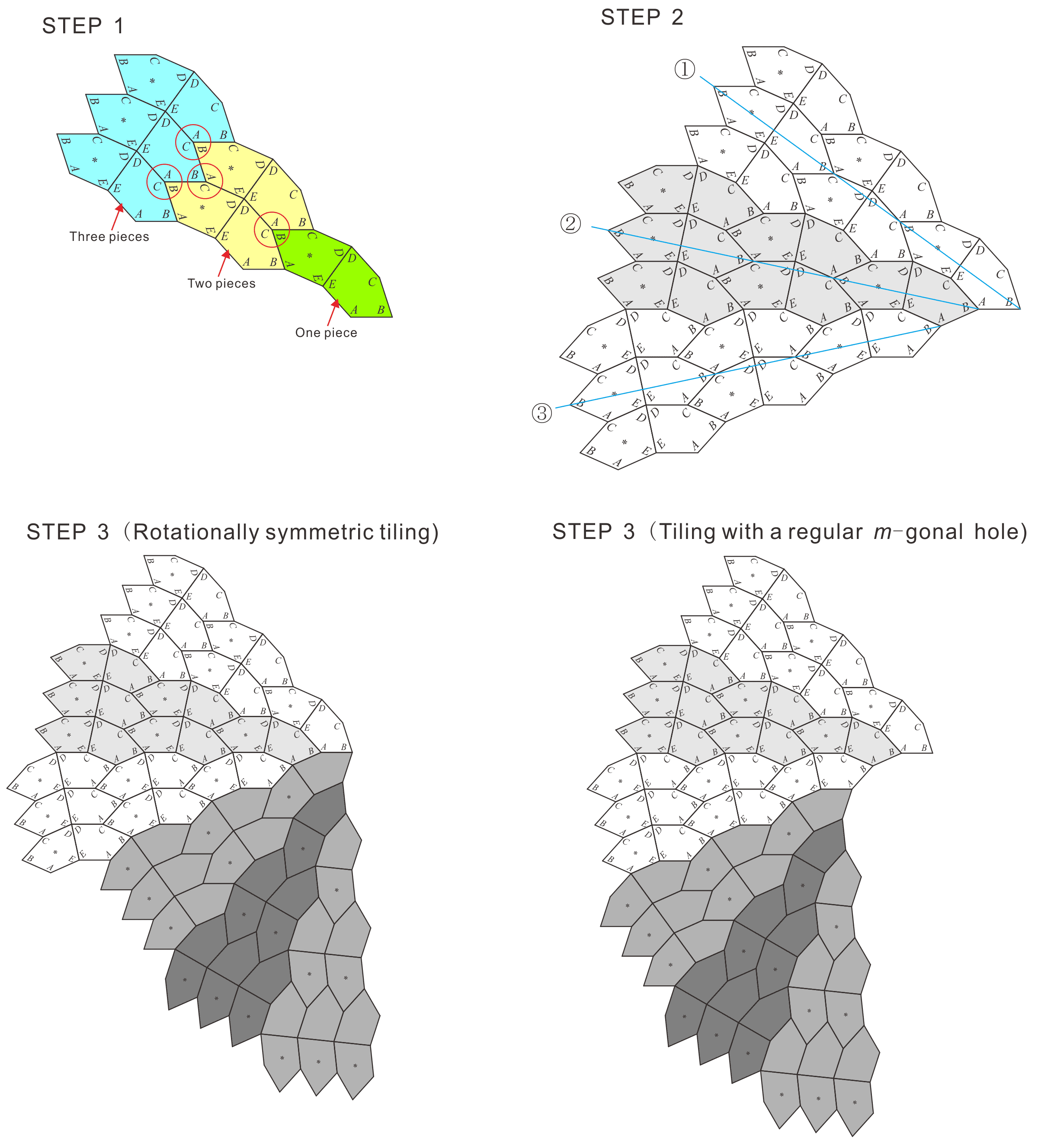} 
  \caption{{\small 
Formation methods of rotationally symmetric edge-to-edge tiling and of 
rotationally symmetric tiling with a regular $m$-gonal hole at the center 
using C11-T1A
} 
\label{fig04}
}
\end{figure}

Here, the formation of rotationally symmetric edge-to-edge tiling with 
convex pentagonal tiles is briefly explained. First, as shown in STEP 1 in 
Figure~\ref{fig04}, create a unit connecting the Octa-units that are generated by 
convex pentagons satisfying (\ref{eq2}) in one direction so that 
$A + B + C = 360^ \circ $. The Octa-units can then be assembled in such a way 
as to increase the number of pieces from one to two to three, and so on, in order. 
Then, copy the two units in STEP 1. In STEP 2, first, rotate one of the two copied 
units by $\frac{120^ \circ }{n}$, with respect to the original (see STEP 2 of 
Figure~\ref{fig04}, wherein the axis of $ \textcircled{\scriptsize 2} $ has a rotation 
of $\frac{120^ \circ }{n}$ with respect to the axis of $ \textcircled{\scriptsize 1} $), 
and connect it to the original to get the following results: $2E + C = 360^ \circ $, 
$2D + A + B = 360^ \circ $, and $A + B + C = 360^ \circ $, as expressed in STEP 2 
of Figure~\ref{fig04}. Next, rotate the remaining unit by $\frac{240^ \circ }{n}$, with 
respect to the original (see STEP 2 of Figure~\ref{fig04}, wherein the axis of 
$ \textcircled{\scriptsize 3} $ has a rotation of $\frac{240^ \circ }{n}$ with respect 
to the axis of $ \textcircled{\scriptsize 1} $), and connect it in the same way as 
the one previously made to make a unit like that in STEP 2 in Figure~\ref{fig04}. 
Next, copy the unit in STEP 2. Subsequently, take the unit from STEP 2 and 
rotate it by the value of the interior angle of vertex $B$. When the original unit 
and the rotated unit are arranged as shown in STEP 3 (Rotationally symmetric 
tiling) in Figure~\ref{fig04}, $\frac{2}{n}$ parts in the $n$-fold rotationally 
symmetric tiling can be formed. Then, by repeating this process as many 
times as necessary, an $n$-fold rotationally symmetric edge-to-edge tiling with 
convex pentagonal tiles can be formed.

Klaassen \cite{Klaassen_2016} presented the convex pentagonal tiles that can 
generate countless rotationally symmetric edge-to-edge tilings; however, the 
pentagon tilings were non-edge-to-edge. On the other hand, it was shown 
that the convex pentagonal tiles satisfying (\ref{eq2}) can generate 
countless rotationally symmetric edge-to-edge tilings. Although the 
corresponding theorem and proof are not presented, similar to how it was 
in \cite{Klaassen_2016}, it will still be understood that countless rotationally 
symmetric tilings can be generated by the foregoing explanations and methods.

%%%%%%%%%%%%%%%%%%%%%%%%%%%%%%%%%%%%%%%%%%%%%%%%%%%%%%%%%%%%%%%%%%%%%%
%%%%%%%%%%%%%%%%%%%%%%%%%%%%%%%%%%%%%%%%%%%%%%%%%%%%%%%%%%%%%%%%%%%%%%
\section{Rotationally symmetric tilings (tiling-like patterns) with a 
regular polygonal hole at the center}
\label{section3}

In \cite{Iliev_2018}, there are figures of rotationally symmetric tiling-like patterns 
with a hole in the center of a regular heptagon, 10-gon, 12-gon, or 18-gon 
formed using convex pentagons (or elements that can be regarded as concave 
octagons). The regular 18-gonal hole can be filled with convex pentagons, 
but the other holes cannot be filled with convex pentagons. 
Note that the tiling-like patterns are not considered tilings due to the 
presence of a gap, but are simply called tilings in this study. The convex 
pentagons forming these rotationally symmetric tilings with a regular 
$m$-gonal hole corresponded to C11-T1A. 

By using the relational expression of the interior angles of (\ref{eq1}), the 
conditions of pentagonal tile C11-T1A that can form rotationally symmetric 
tilings with a regular $m$-gonal hole are expressed in (\ref{eq3}).

\begin{equation}
\label{eq3}
\left\{ {\begin{array}{l}
 A = 180^ \circ - \dfrac{720^ \circ}{m\strut}, \\ 
 B = \dfrac{1080^ \circ}{m\strut}, \\ 
 C = 180^ \circ - \dfrac{360^ \circ}{m\strut}, \\ 
 D = 90^ \circ - \dfrac{180^ \circ}{m\strut}, \\ 
 E = 90^ \circ + \dfrac{180^ \circ}{m\strut}, \\ 
 a = b = c = d, \\ 
 \end{array}} \right.
\end{equation}

\noindent
where $m$ is an integer greater than or equal to seven, because 
$0^ \circ < B < 180^ \circ $. Note that, because ``$180^ \circ - \frac{360^ \circ}{m}$" 
corresponds to one interior angle of a regular $m$-gon, the value of 
``$A+B$" in (\ref{eq3}) is equal to the outer angle ($180^ \circ + \frac{360^ \circ}{m}$) 
of one vertex of a regular $m$-gon.

Table~\ref{tab2} presents some of the relationships between the interior angles of 
convex pentagons satisfying (\ref{eq3}) that can form the rotationally symmetric 
tilings with a regular $m$-gonal hole at the center. (For $m = 7\!-\!10, 12, 14, 15, 
18, 21, 24, 27$, tilings with a regular $m$-gonal hole at  the center formed by 
convex pentagons are drawn. For further details, Figures~\ref{fig18}--\ref{fig26}.) 
If these elements are considered to be convex pentagons, the 
connection is edge-to-edge. These tilings with a regular $m$-gonal hole with 
$D_{m}$ symmetry\footnote{ 
``$D_{m}$" is based on the Schoenflies notation for symmetry in a two-dimensional 
point group~\cite{wiki_point_group, wiki_schoenflies_notation}. ``$D_{m}$" represents 
an $m$-fold rotation axis with $m$ reflection symmetry axes.
} at the center have $C_{m}$ symmetry. If convex pentagons satisfying (\ref{eq3}) 
have $m$ that is divisible by three, they are also convex pentagonal tiles that 
satisfy (\ref{eq2}).

\begin{table}[htb]
 \begin{center}
{\small
\caption[Table 2]{
Interior angles of convex pentagons satisfying (\ref{eq3}) that can form the 
rotationally symmetric tilings with a regular $m$-gonal hole at the center
}
\label{tab2}
}
\begin{tabular}
{c| D{.}{.}{2} D{.}{.}{2} D{.}{.}{2} D{.}{.}{2} D{.}{.}{2} |c|c|c}
%{c|rrrrr|c|c}
\hline
\raisebox{-1.50ex}[0cm][0cm]{$n$}& 
\multicolumn{5}{c|}{\shortstack{ Value of interior angle (degree) } } & 
\raisebox{-4.6ex}[0.7cm][0.5cm]{\shortstack{Edge \\length \\of $e$}} & 
\raisebox{-3.0ex}[0.7cm][0.5cm]{\shortstack{$n$ of \\Table~\ref{tab1}}} & 
\raisebox{-3.0ex}[0.7cm][0.5cm]{\shortstack{Figure \\number}} \\

 & 
\textit{A} & 
\textit{B}& 
\textit{C}& 
\textit{D}& 
\textit{E}& 
 & 
  \\
\hline
7& 
77.14 & 
154.29 & 
128.57 & 
64.29 & 
115.71 & 
1.950 & 
& 
\ref{fig16} \\
\hline
8& 
90& 
135& 
135& 
67.5& 
112.5& 
1.848 & 
& 
\ref{fig17} \\
\hline
9& 
100& 
120& 
140& 
70& 
110& 
1.732 & 
3& 
\ref{fig18} \\
\hline
10& 
108& 
108& 
144& 
72& 
108& 
1.618 & 
& 
\ref{fig19} \\
\hline
11& 
114.55 & 
98.18 & 
147.27 & 
73.64 & 
106.36 & 
1.511 & 
& 
 \\
\hline
12& 
120& 
90& 
150& 
75& 
105& 
1.414 & 
4& 
\ref{fig20} \\
\hline
13& 
124.62 & 
83.08 & 
152.31 & 
76.15 & 
103.85 & 
1.326 & 
& 
 \\
\hline
14& 
128.57 & 
77.14 & 
154.29 & 
77.14 & 
102.86 & 
1.247 & 
& 
\ref{fig21} \\
\hline
15& 
132& 
72& 
156& 
78& 
102& 
1.176 & 
5& 
\ref{fig22} \\
\hline
16& 
135& 
67.5& 
157.5& 
78.75& 
101.25& 
1.111 & 
& 
 \\
\hline
17& 
137.65 & 
63.53 & 
158.82 & 
79.41 & 
100.59 & 
1.053 & 
& 
 \\
\hline
18& 
140& 
60& 
160& 
80& 
100& 
1 & 
6& 
\ref{fig23} \\
\hline
19& 
142.11 & 
56.84 & 
161.05 & 
80.53 & 
99.47 & 
0.952 & 
& 
 \\
\hline
20& 
144& 
54& 
162& 
81& 
99& 
0.908 & 
& 
 \\
\hline
21& 
145.71 & 
51.43 & 
162.86 & 
81.43 & 
98.57 & 
0.868 & 
7& 
\ref{fig25} \\
\hline
22& 
147.27 & 
49.09 & 
163.64 & 
81.82 & 
98.18 & 
0.831 & 
& 
 \\
\hline
23& 
148.70 & 
46.96 & 
164.35 & 
82.17 & 
97.83 & 
0.797 & 
& 
 \\
\hline
24& 
150& 
45& 
165& 
82.5& 
97.5& 
0.765 & 
8& 
\ref{fig26} \\
\hline
25& 
151.2& 
43.2& 
165.6& 
82.8& 
97.2& 
0.736 & 
& 
 \\
\hline
26& 
152.31 & 
41.54 & 
166.15 & 
83.08 & 
96.92 & 
0.709 & 
& 
 \\
\hline
27& 
153.33 & 
40 & 
166.67 & 
83.33 & 
96.67 & 
0.684 & 
9& 
\ref{fig27}
 \\
\hline
...& 
...& 
...& 
...& 
...& 
...& 
...& 
 \\
\hline
\end{tabular}

\end{center}

\end{table}

Rotationally symmetric tilings with a regular $m$-gonal hole at the center can 
be made in almost the same way as the $n$-fold rotationally symmetric 
edge-to-edge tiling (i.e., the manner of creation follows the same steps up 
to STEP 2 of Figure~\ref{fig04}). The part with three edges \textit{AB}, at the 
center of the unit made in STEP 2 of Figure~\ref{fig04}, corresponds to the 
contour of a regular $m$-gon. Therefore, copy the unit in STEP 2, take the 
unit from STEP 2 and rotate it by the value of the interior angle of vertex $B$. 
When the original unit and the rotated unit are arranged as shown in STEP 3 
(Tiling with a regular $m$-gonal hole) in Figure~\ref{fig04}, $\frac{6}{m}$ parts of a 
regular $m$-gon can be formed. Then, by repeating this process as many 
times as necessary, a rotationally symmetric tiling with a regular $m$-gonal 
hole at the center can be formed.

As shown in Table~\ref{tab2} and Figure~\ref{fig23}, the convex pentagonal tile that 
satisfies (\ref{eq3}), where $m = 18$, is an equilateral convex pentagon. As such, 
the regular 18-gonal hole can be filled with equilateral convex pentagons. The 
pentagonal arrangement pattern in the regular 18-gon is unique; however, the 
regular 18-gon can be reversed; therefore, two patterns can be used in the 
hole (see Figure~\ref{fig23}). The tiling in Figure~\ref{fig24} is one of the six-fold 
rotationally symmetric edge-to-edge tilings with equilateral convex pentagons 
that Hirschhorn and Hunt present in \cite{H_and_H_1985}. In addition, the tiling in 
Figure~\ref{fig23} is that of Hirschhorn~\cite{Klaassen_2017}. In \cite{Iliev_2018}, there 
are figures of tilings corresponding to Figures~\ref{fig08}, \ref{fig23}, and \ref{fig24}, 
with the one corresponding to Figure~\ref{fig24} labeled as Hirschhorn's rosette. 
These tilings with a regular 18-gonal hole that is filled with equilateral convex 
pentagons have $C_{6}$ symmetry. In Pegg's research~\cite{Pegg_site}, they 
stated there exists the tiling of Figure~\ref{fig08}, using the equilateral convex 
pentagon, and that there are also other interesting tilings, such as that 
with the ``crystal-like'' structure. Note that in Pegg's research~\cite{Pegg_site}, 
there are non-monohedral tilings created using convex and concave pentagons.

As described above, the tiling of $m = 10$ in Table~\ref{tab2} (see Figure~\ref{fig19}) is 
listed in Iliev's research \cite{Iliev_2018}, however, the value of the interior angle of 
the convex pentagon was erroneous, as well as in Smith's research \cite{Smith_site}. 
It is shown that the convex pentagon of $m = 10$ in Table~\ref{tab2} is the same as the 
convex pentagonal tile belonging to both the Type 1 and Type 2 families, and 
that it can form various tilings~\cite{Bailey_site, Iliev_2018, Smith_site}. 
In \cite{Bailey_site} and \cite{Smith_site}, there are also 
tilings with a regular pentagonal hole at the center using this convex 
pentagon of $m = 10$ in Table~\ref{tab2}. Moreover, it is interesting to note that 
the ring-shaped layer forming the periphery can be reversed in \cite{Smith_site}.

%%%%%%%%%%%%%%%%%%%%%%%%%%%%%%%%%%%%%%%%%%%%%%%%%%%%%%%%%%%%%%%%%%%%%%
%%%%%%%%%%%%%%%%%%%%%%%%%%%%%%%%%%%%%%%%%%%%%%%%%%%%%%%%%%%%%%%%%%%%%%
\section{Spiral tilings}
\label{section4}

Klaassen \cite{Klaassen_2016} also showed the viewpoint of rotationally symmetric 
tilings with convex pentagons as a spiral structure. A similar spiral structure can 
be found in the rotationally symmetric edge-to-edge tilings with convex 
pentagons of C11-T1A. Figure~\ref{fig28} pertains to the case where $n = 5$, 
which makes the spiral structure easier to understand. From this spiral structure, 
the difference between the rotationally symmetric tilings with convex 
pentagons in \cite{Klaassen_2016} and the rotationally symmetric edge-to-edge 
tilings with C11-T1A other than the edge-to-edge property will be determinable.

Here, let us pay attention to the convex pentagon of $m = 8$ in Table~\ref{tab2}. 
This convex pentagon has the property of $B = C = 135^ \circ $. From this 
property, the convex pentagon of $m = 8$ in Table~\ref{tab2} can form spiral 
tilings\footnote{ 
It is difficult to determine what the spiral tiling is. There are some discussions, 
such as those in \cite{Klaassen_2017}. The tiling shown in Figure~\ref{fig28} can 
be regarded as a spiral, however, it has a radial structure extending from the 
center. Thus, we feel that it is more straightforward to regard it as a rotationally 
symmetric tiling. On the other hand, because the tiling of Figure~\ref{fig29} 
does not have a radial structure extending from the center or singularity point, 
we feel that it is a spiral rather than a rotationally symmetric tiling.
} with two-fold rotational symmetry, as shown in Figure~\ref{fig29}(a). 
This convex pentagonal tile belongs to both the Type 1 and Type 7 
families~\cite{G_and_S_1987, Sugimoto_NoteTP, Sugimoto_2012, 
Sugimoto_2016, wiki_pentagon_tiling}.

After finding the tiling in Figure~\ref{fig29}(a), we found that, in ``(6) Central 
block `wraparound'~'' of \cite{Bailey_site}, Bailey presented the spiral tilings with two-fold 
rotational symmetry by the convex pentagon corresponding to the case of $m = 
10$ in Table~\ref{tab2}. The spiral tiling is the same as that in Figure~\ref{fig30}(a). 
Furthermore, from \cite{Bailey_site}, we found that it is possible to form tilings that 
keep the spiral structure and extend in one direction only, as shown in 
Figures~\ref{fig30}(b), \ref{fig30}(c), and \ref{fig30}(d). Then, we determined that the 
convex pentagon of $m = 8$ in Table~\ref{tab2} can also form tilings that keep the spiral 
structure and extend in one direction (see Figures~\ref{fig29}(b), \ref{fig29}(c), and 
\ref{fig29}(d)).

The convex pentagon of $m = 14$ in Table~\ref{tab2} has $B = D = \frac{3\pi }{7} 
\approx 77.14^ \circ $, and belongs to both the Type 1 and Type 2 families. 
Similar to the convex pentagons of $m = 8, 10$ in Table~\ref{tab2}, the convex 
pentagon of $m = 14$ in Table~\ref{tab2} can also form the spiral tiling with 
two-fold rotational symmetry (see Figure~\ref{fig31}(a)) and tilings that keep the 
spiral structure and extend in one direction (see Figures~\ref{fig31}(b), \ref{fig31}(c), 
and \ref{fig31}(d)). For all convex pentagons of $m = 8, 10, 14$ in Table~\ref{tab2}, 
it is also possible to remove one spiral structure and to extend the belts wherein the 
Octa-units are arranged (see Figure~\ref{fig32}).

It should be noted that the spiral tiling with a regular pentagonal hole at 
the center is shown in \cite{Smith_site} using the convex pentagon corresponding 
to $m = 10$ in Table~\ref{tab2}.

%%%%%%%%%%%%%%%%%%%%%%%%%%%%%%%%%%%%%%%%%%%%%%%%%%%%%%%%%%%%%%%%%%%%%%
%%%%%%%%%%%%%%%%%%%%%%%%%%%%%%%%%%%%%%%%%%%%%%%%%%%%%%%%%%%%%%%%%%%%%%
\section{Conclusions}

The existence of convex pentagonal tiles that can generate countless number of 
rotationally symmetric edge-to-edge tilings is not well-known. As mentioned 
above, Klaassen presented a case where tilings are not edge-to-edge. From 
\cite{Iliev_2018}, we determined that the convex pentagon called C11-T1A 
in \cite{S_and_O_2009} had this property. In \cite{Iliev_2018}, Iliev presented 
the tilings corresponding to $n = 4, 5, 6, 10$ in Table~\ref{tab1} and 
$m = 7, 10, 12, 18$ in Table~\ref{tab2}. It may have been inferred by 
Iliev that there are countless rotationally symmetric tilings similar 
to those; however, there seems to be no clarification on the matter. 
In particular, there is no description of the concrete properties 
(tile conditions), such as the convex pentagon of C11-T1A, as shown in 
this study.

It was previously recognized that the convex pentagon belonging to both the 
Type 1 and Type 7 families (i.e., convex pentagon corresponding to $m = 8$ in 
Table~\ref{tab2}) can generate various albeit limited tilings. However, in reality, 
we noticed that it is possible to generate more tilings. Among them are 
other rotational symmetric tilings and spiral tilings. We will introduce 
them in a different article~\cite{Sugimoto_2020_4} .

\bigskip
\noindent
\textbf{Acknowledgments.} 
The author would like to thank Yoshiaki ARAKI of Japan Tessellation Design 
Association, for information of Bailey's site, etc. and comments.

%%%%%%%%%%%%%%%%%%%%%%%%%%%%%%%%%%%%%%%%%%%%%%%%%%%%%%%%%%%%%%%%%%%%%%
%%%%%%%%%%%%%%%%%%%%%%%%%%%%%%%%%%%%%%%%%%%%%%%%%%%%%%%%%%%%%%%%%%%%%%

\vspace{5cm}

%%%%%%%%%%%%%%%%%%%%%%%%%%%%%%%%%%%%%%%%%%%%%%%%%%%%%%%%%%%%%%%%%%%%%%
%%%%%%%%%%%%%%%%%%%%%%%%%%%%%%%%%%%%%%%%%%%%%%%%%%%%%%%%%%%%%%%%%%%%%%

\renewcommand{\figurename}{{\small Figure.}}
\begin{figure}[htbp]
 \centering\includegraphics[width=10cm,clip]{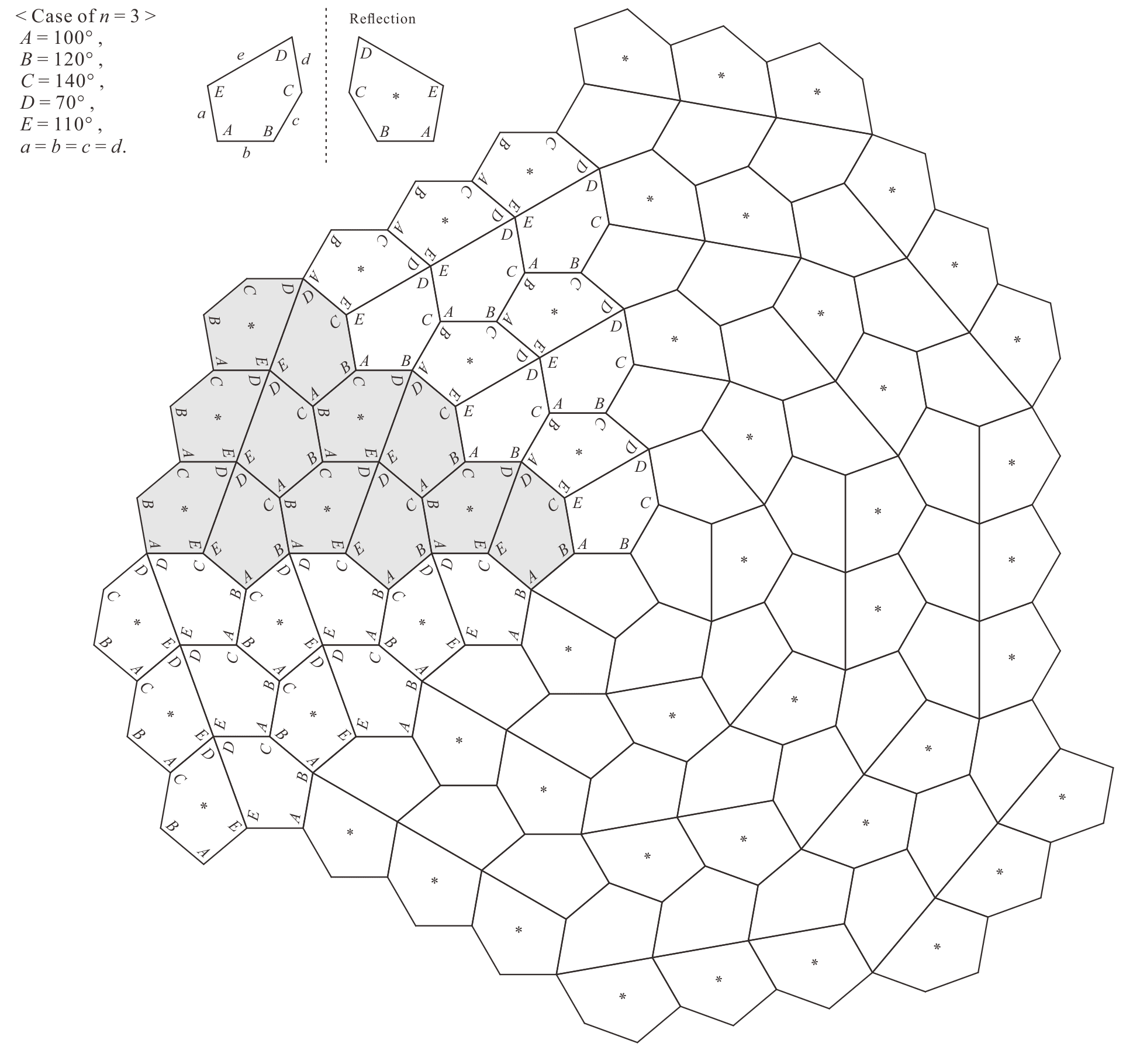} 
  \caption{{\small 
Three-fold rotationally symmetric edge-to-edge tiling by a convex 
pentagon 
} 
\label{fig05}
}
\end{figure}

\renewcommand{\figurename}{{\small Figure.}}
\begin{figure}[htbp]
 \centering\includegraphics[width=10cm,clip]{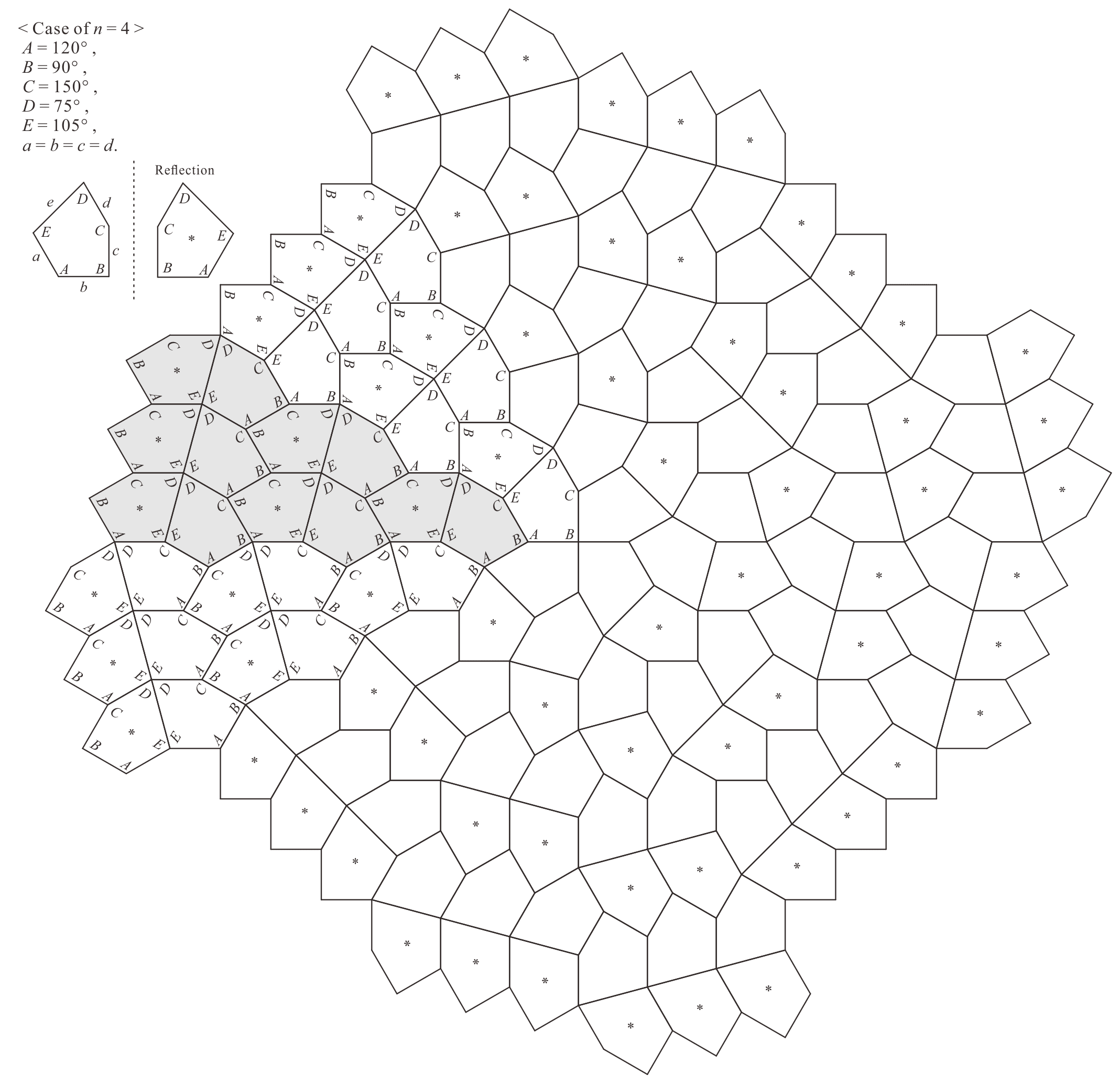} 
  \caption{{\small 
Four-fold rotationally symmetric edge-to-edge tiling by a convex 
pentagon 
} 
\label{fig06}
}
\end{figure}

\renewcommand{\figurename}{{\small Figure.}}
\begin{figure}[htbp]
 \centering\includegraphics[width=10cm,clip]{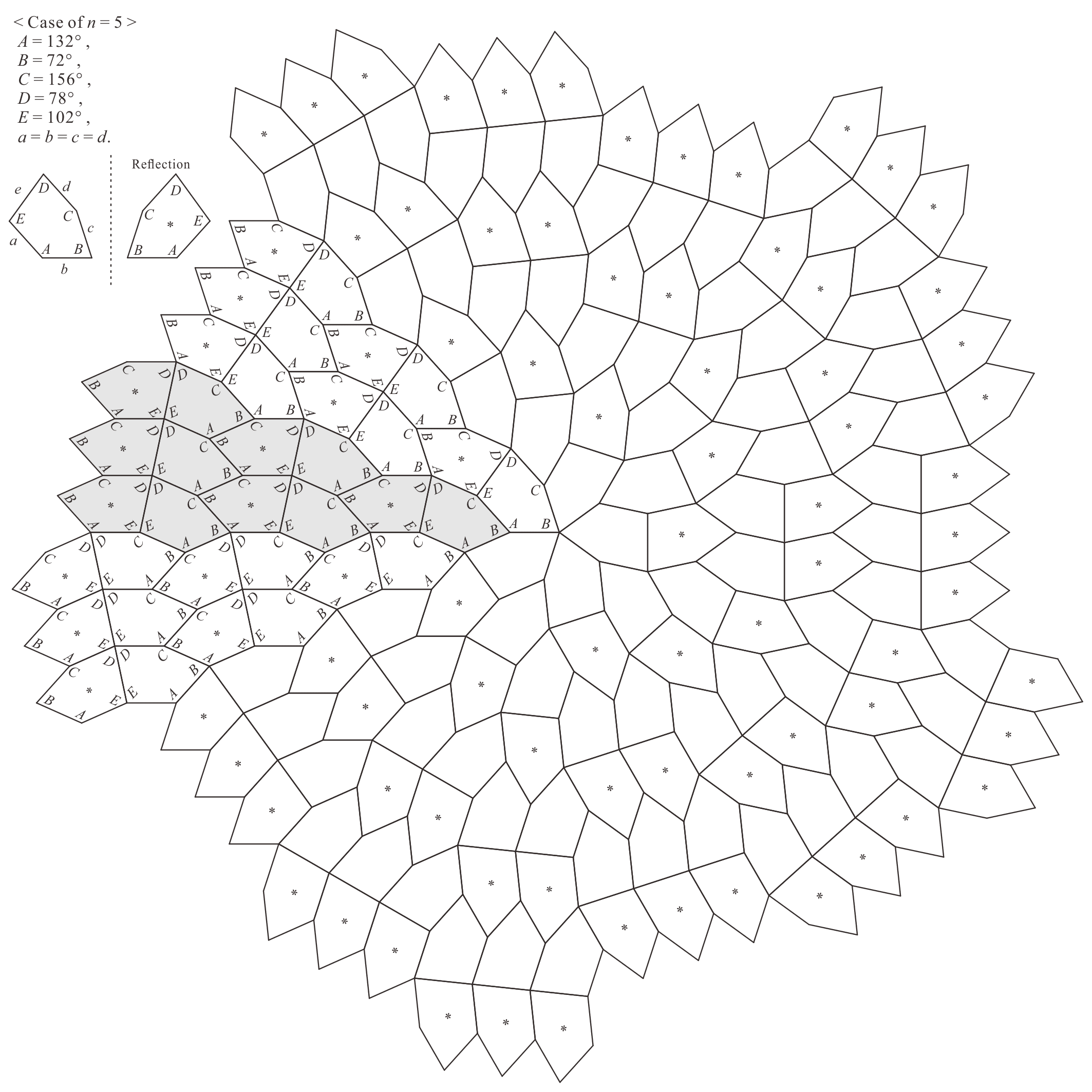} 
  \caption{{\small 
Five-fold rotationally symmetric edge-to-edge tiling by a convex 
pentagon 
} 
\label{fig07}
}
\end{figure}

\renewcommand{\figurename}{{\small Figure.}}
\begin{figure}[htbp]
 \centering\includegraphics[width=10cm,clip]{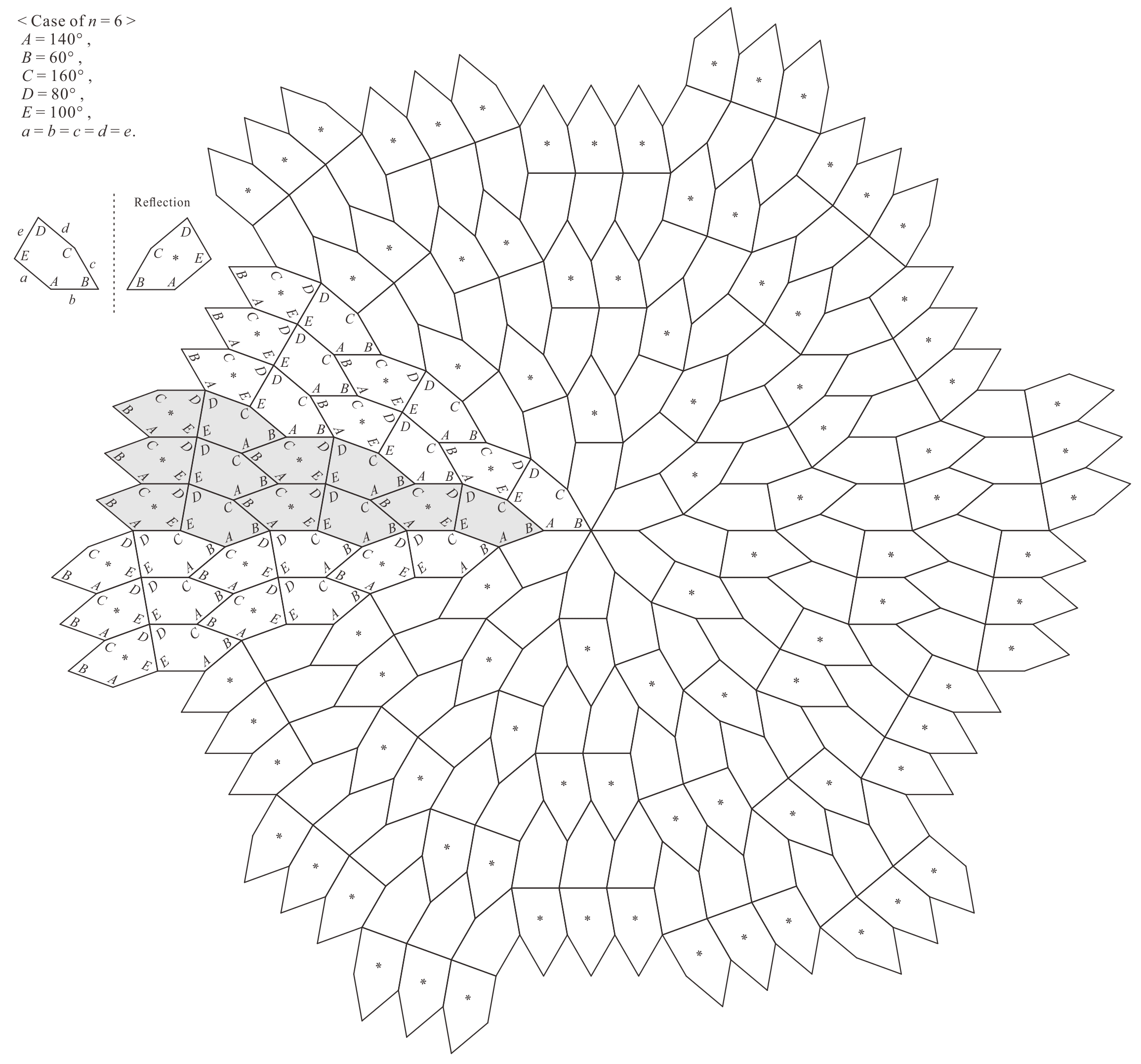} 
  \caption{{\small 
Six-fold rotationally symmetric edge-to-edge tiling by a convex 
pentagon 
} 
\label{fig08}
}
\end{figure}

\renewcommand{\figurename}{{\small Figure.}}
\begin{figure}[htbp]
 \centering\includegraphics[width=10cm,clip]{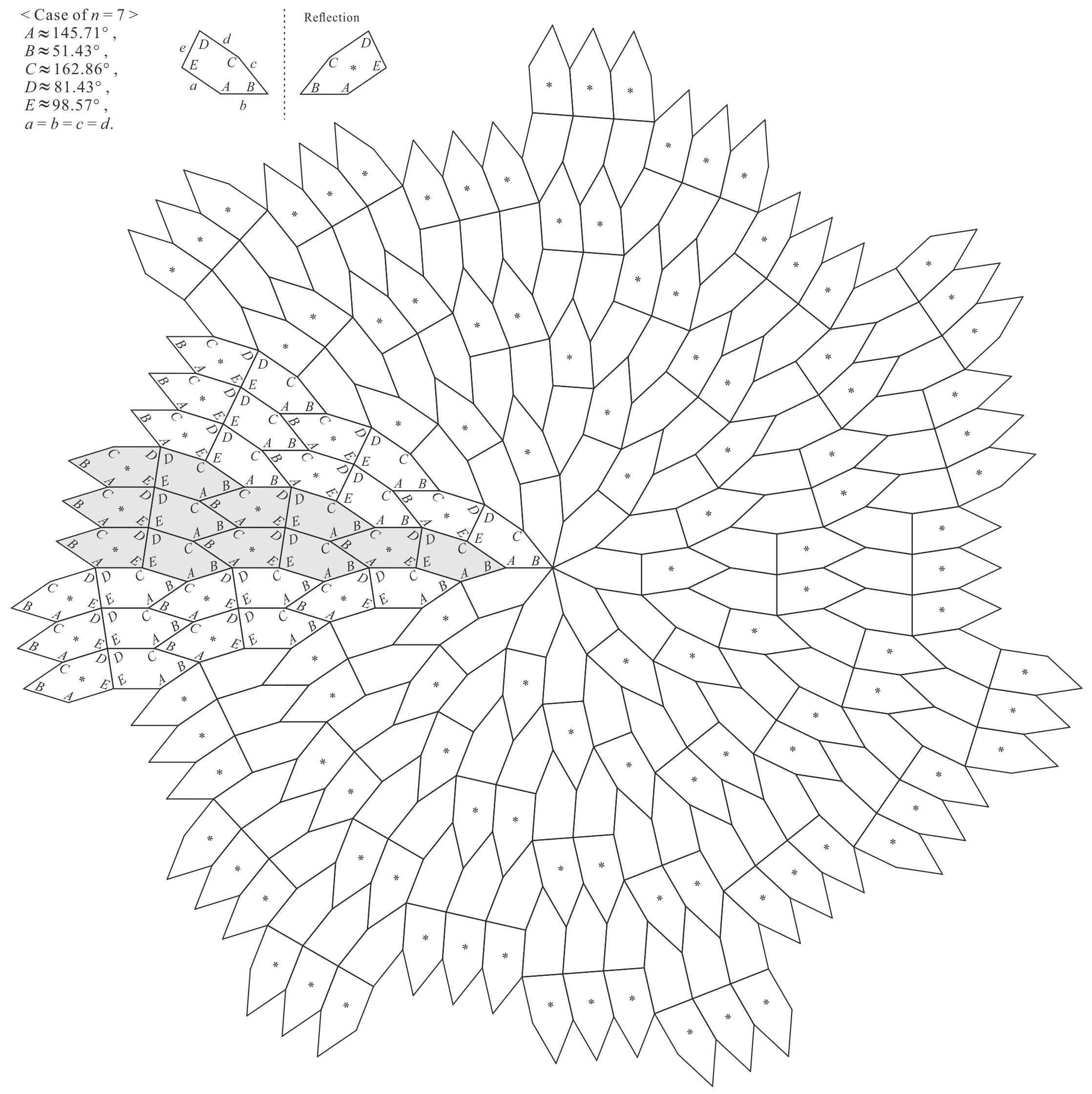} 
  \caption{{\small 
Seven-fold rotationally symmetric edge-to-edge tiling by a convex 
pentagon 
} 
\label{fig09}
}
\end{figure}

\renewcommand{\figurename}{{\small Figure.}}
\begin{figure}[htbp]
 \centering\includegraphics[width=10cm,clip]{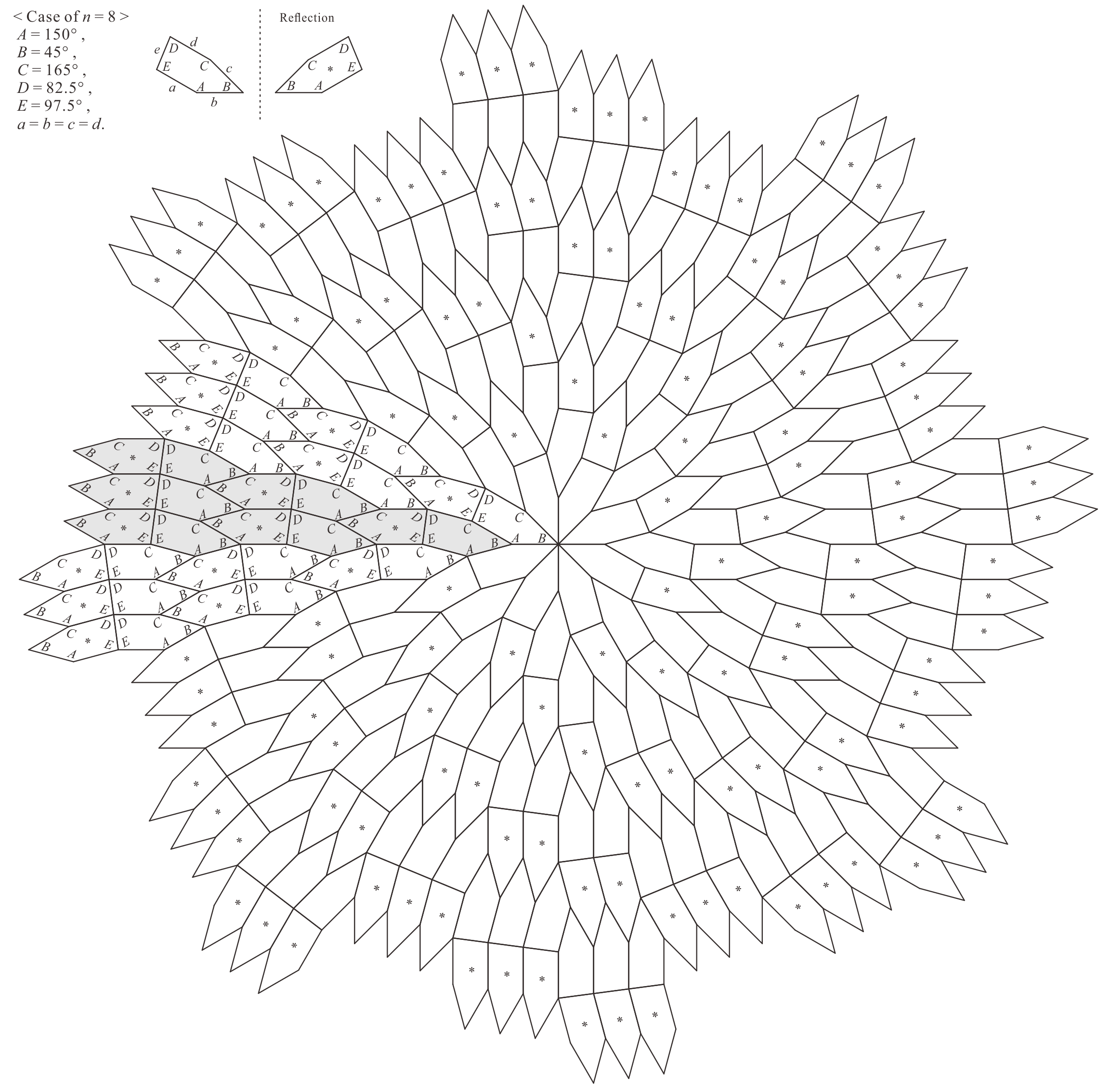} 
  \caption{{\small 
Eight-fold rotationally symmetric edge-to-edge tiling by a convex 
pentagon 
} 
\label{fig10}
}
\end{figure}

\renewcommand{\figurename}{{\small Figure.}}
\begin{figure}[htbp]
 \centering\includegraphics[width=10cm,clip]{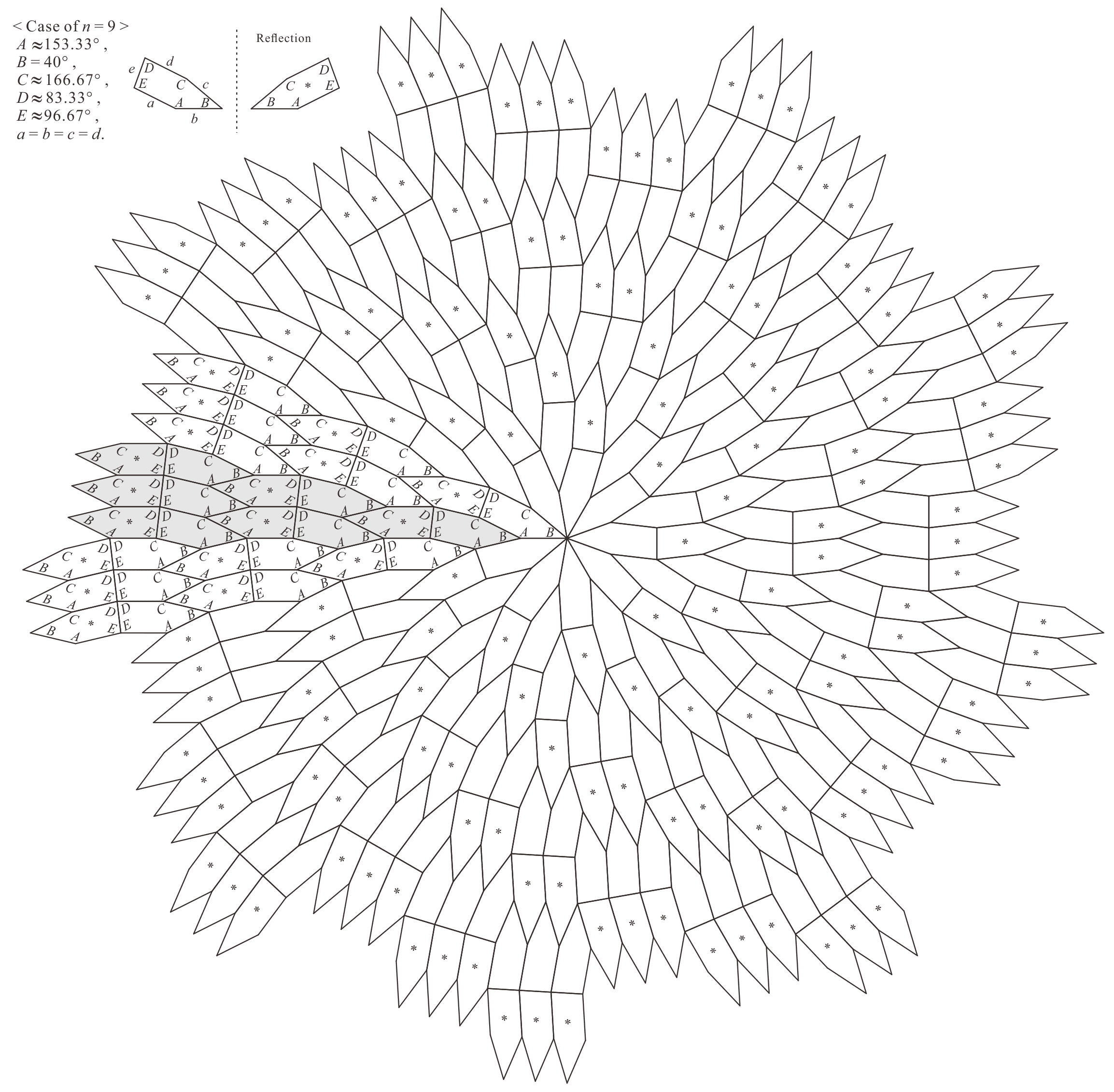} 
  \caption{{\small 
Nine-fold rotationally symmetric edge-to-edge tiling by a convex 
pentagon 
} 
\label{fig11}
}
\end{figure}

\renewcommand{\figurename}{{\small Figure.}}
\begin{figure}[htbp]
 \centering\includegraphics[width=10cm,clip]{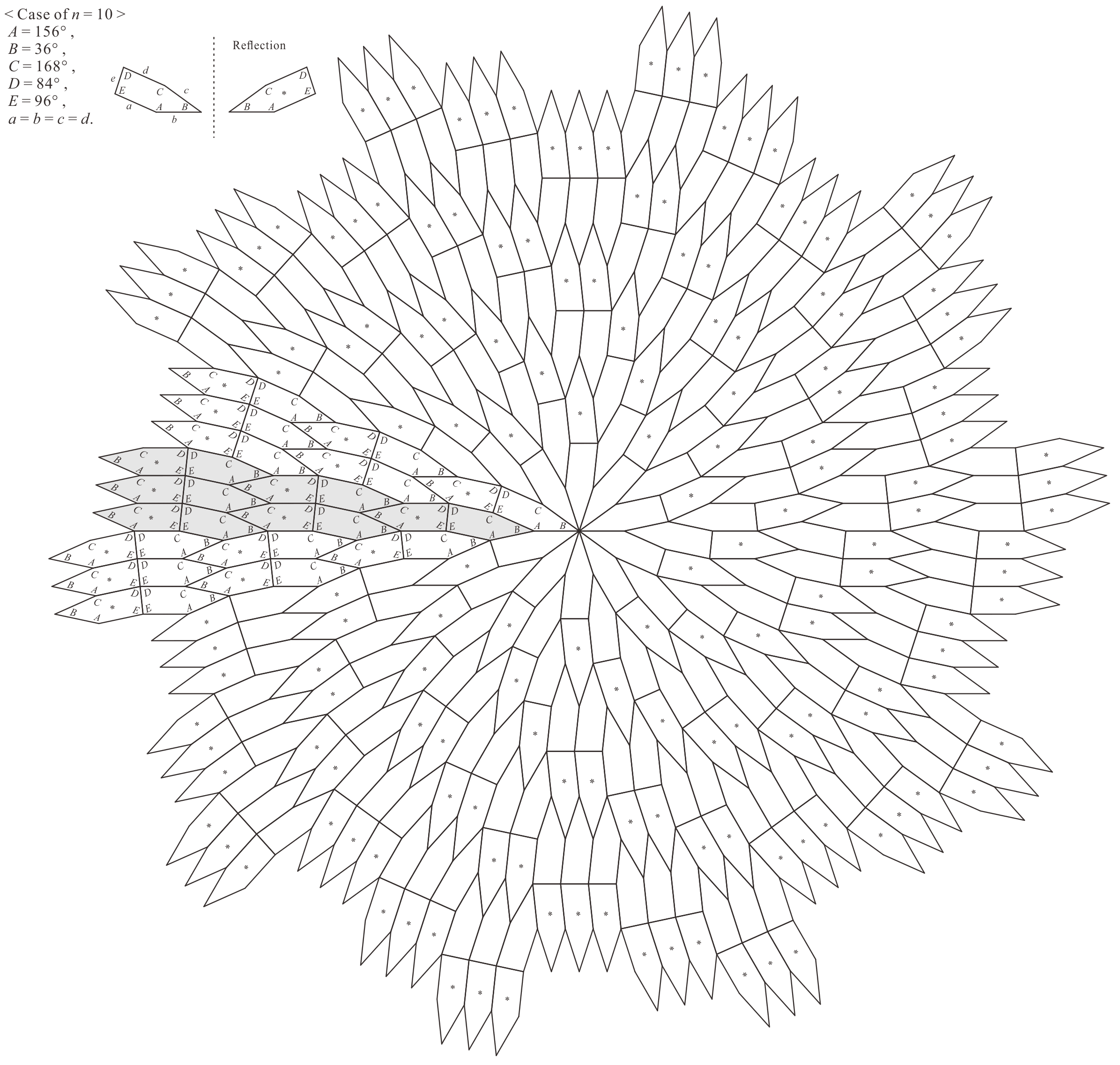} 
  \caption{{\small 
10-fold rotationally symmetric edge-to-edge tiling by a convex 
pentagon 
} 
\label{fig12}
}
\end{figure}

\renewcommand{\figurename}{{\small Figure.}}
\begin{figure}[htbp]
 \centering\includegraphics[width=10cm,clip]{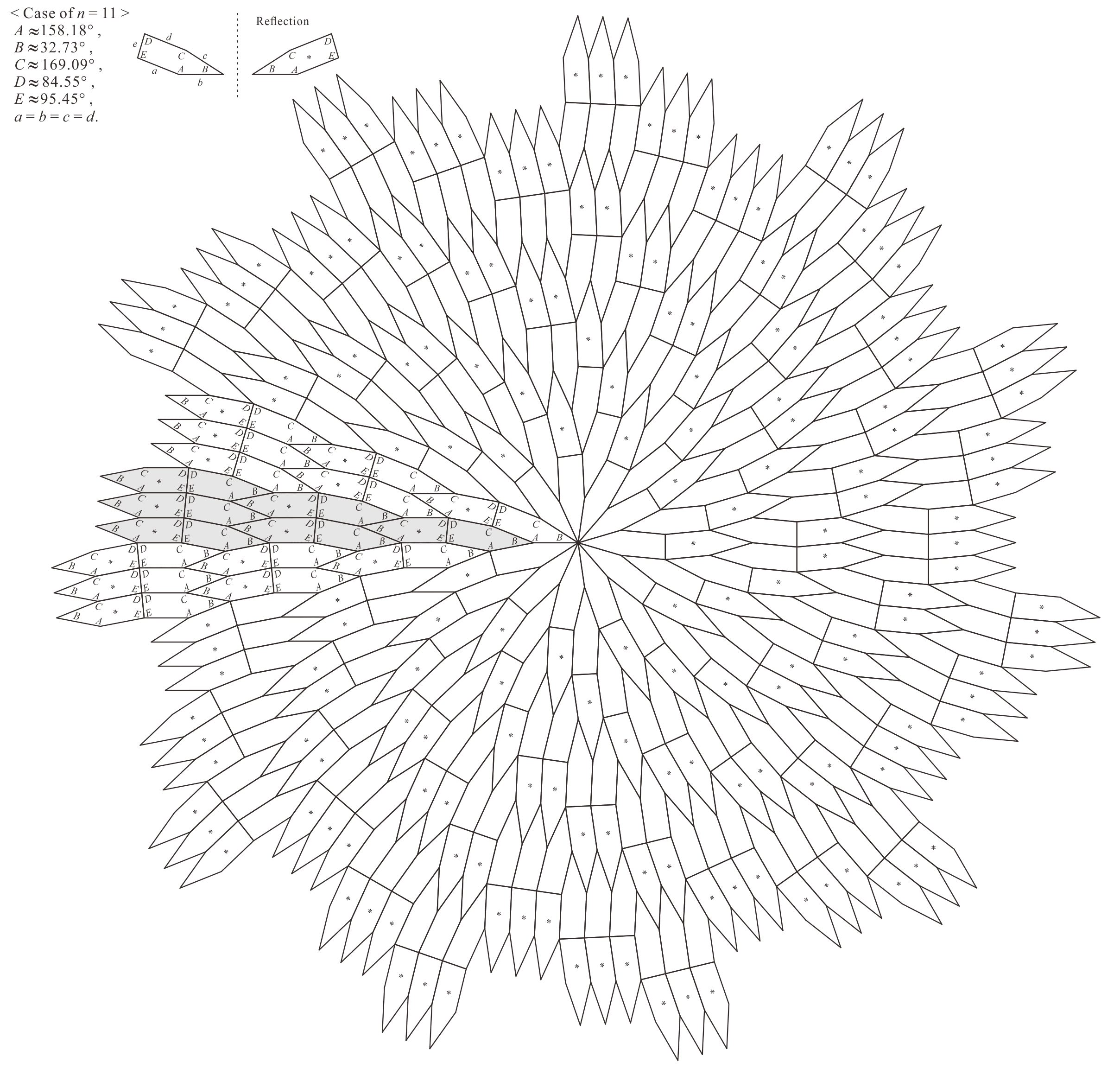} 
  \caption{{\small 
11-fold rotationally symmetric edge-to-edge tiling by a convex 
pentagon 
} 
\label{fig13}
}
\end{figure}

\renewcommand{\figurename}{{\small Figure.}}
\begin{figure}[htbp]
 \centering\includegraphics[width=10cm,clip]{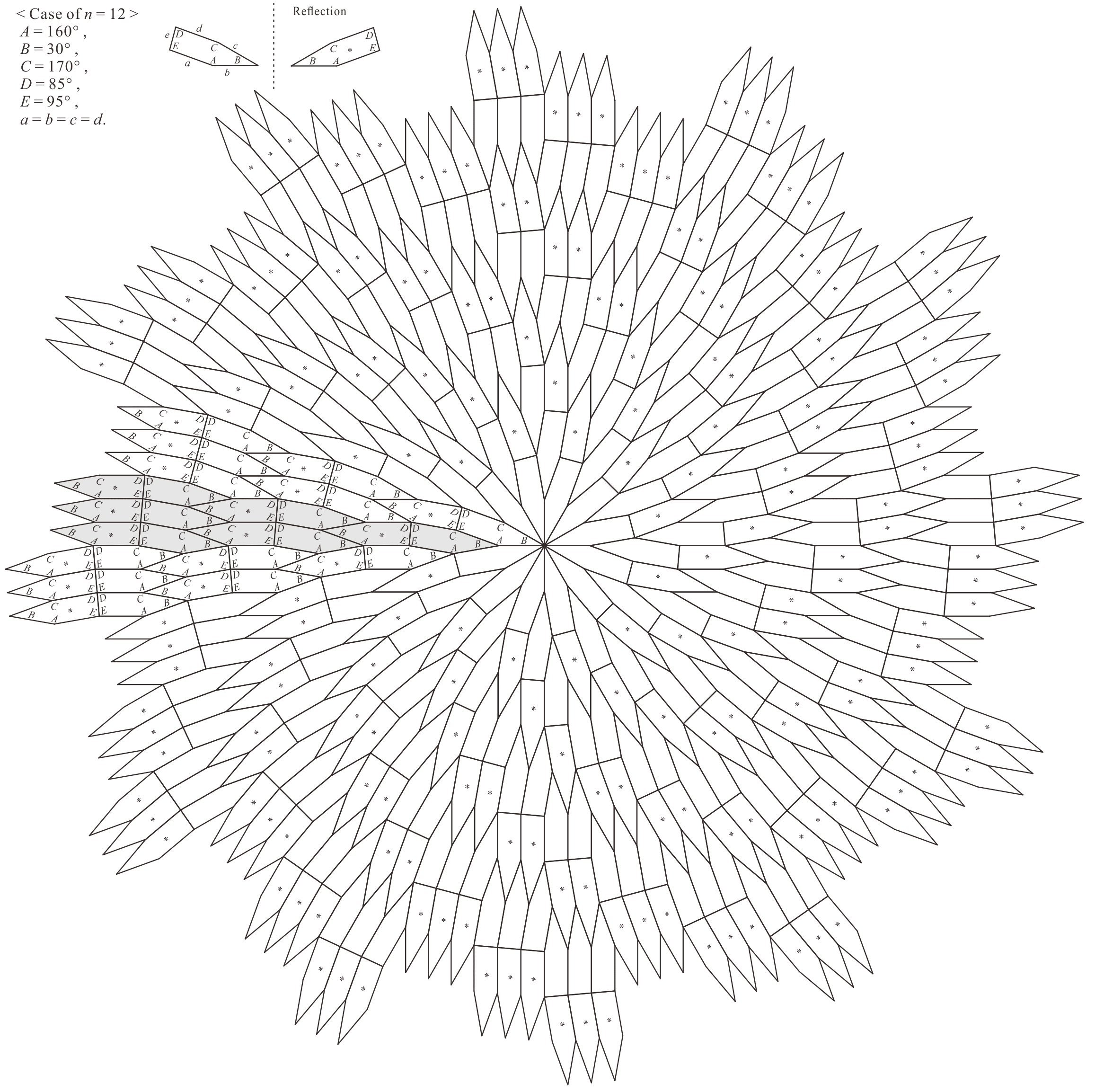} 
  \caption{{\small 
12-fold rotationally symmetric edge-to-edge tiling by a convex 
pentagon 
} 
\label{fig14}
}
\end{figure}

\renewcommand{\figurename}{{\small Figure.}}
\begin{figure}[htbp]
 \centering\includegraphics[width=10cm,clip]{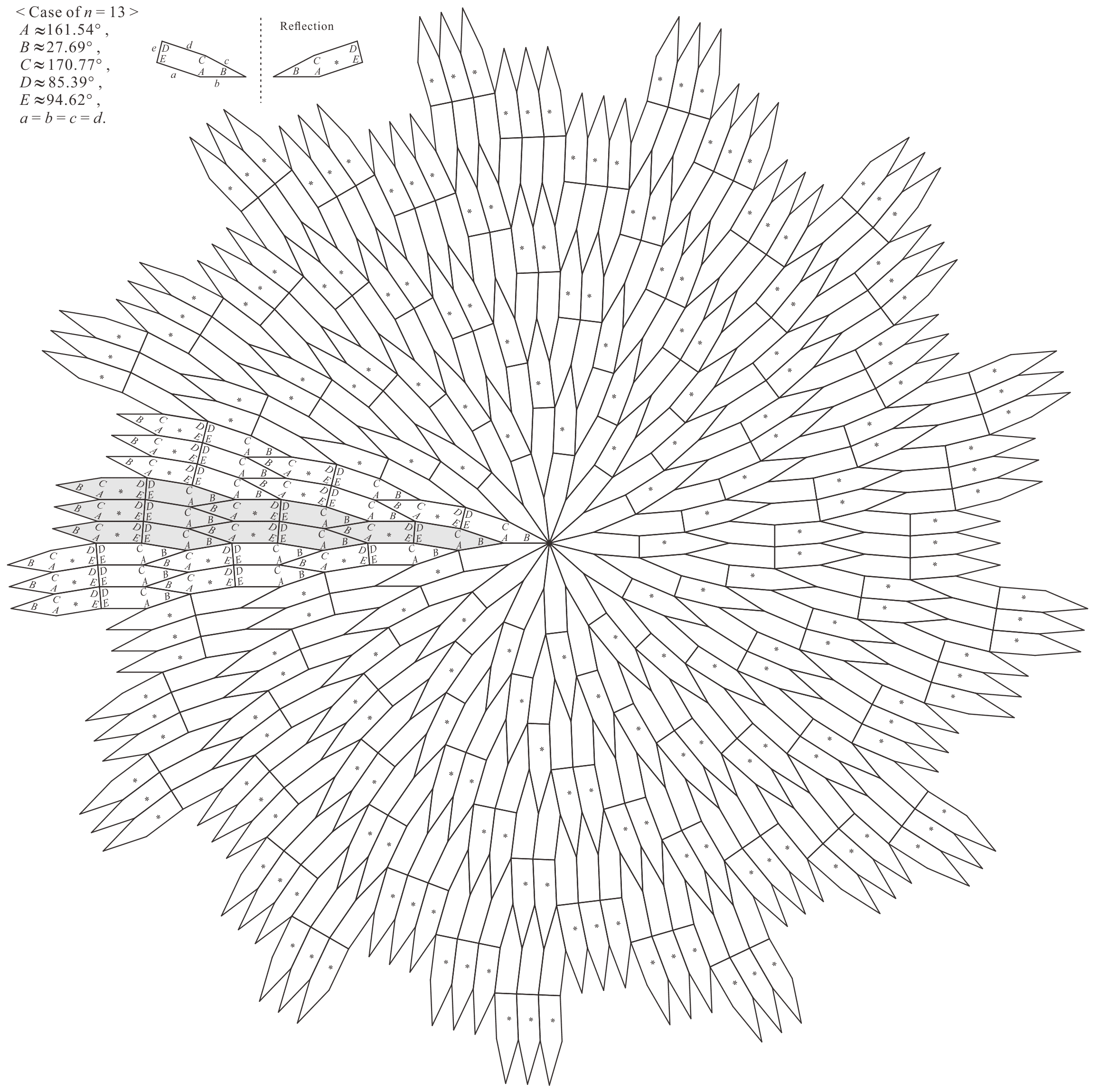} 
  \caption{{\small 
13-fold rotationally symmetric edge-to-edge tiling by a convex 
pentagon
} 
\label{fig15}
}
\end{figure}

\renewcommand{\figurename}{{\small Figure.}}
\begin{figure}[htbp]
 \centering\includegraphics[width=11cm,clip]{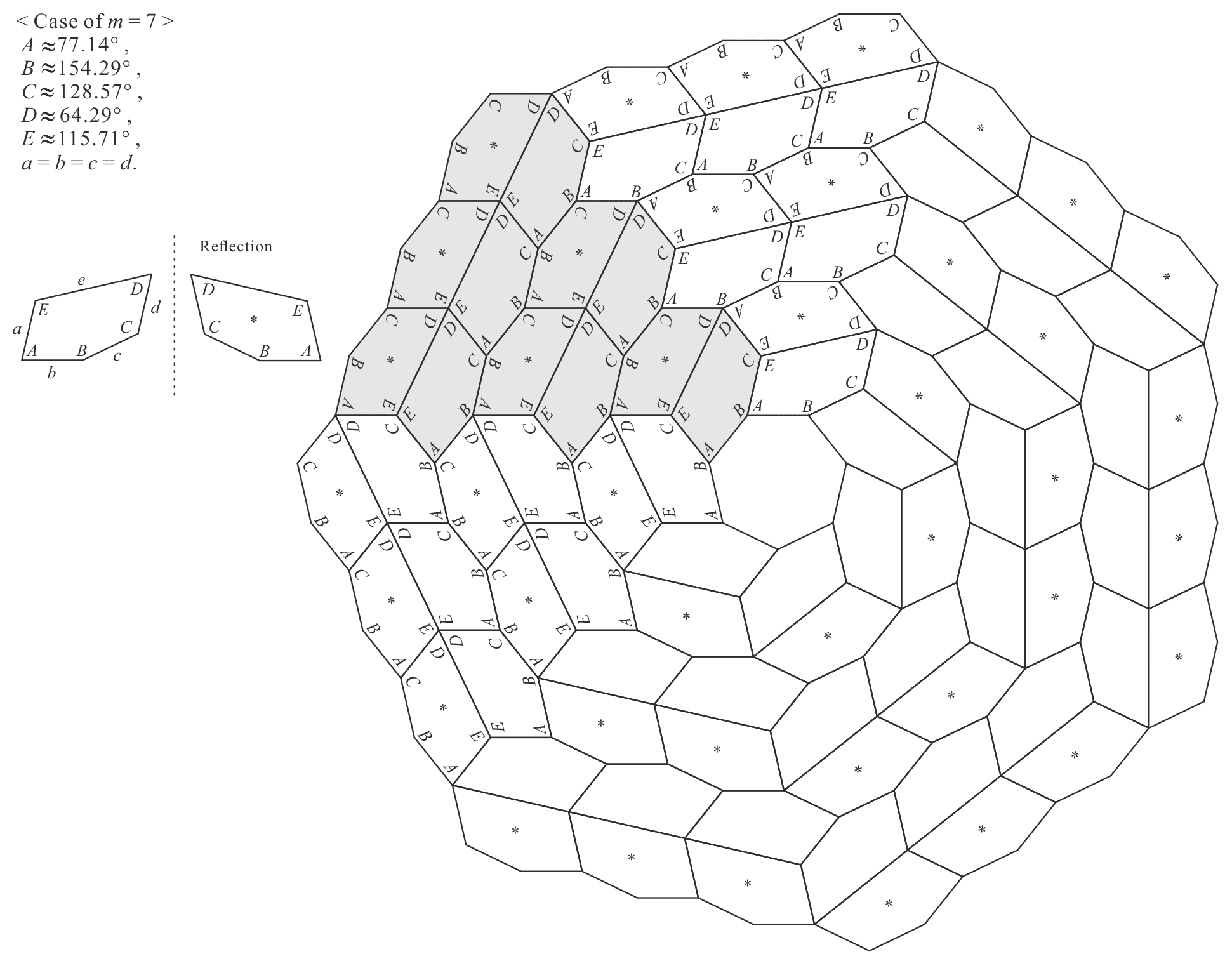} 
  \caption{{\small 
Rotationally symmetric tiling with $C_{7}$ symmetry with a regular 
convex heptagonal hole at the center by a convex pentagon
} 
\label{fig16}
}
\end{figure}

\renewcommand{\figurename}{{\small Figure.}}
\begin{figure}[htbp]
 \centering\includegraphics[width=11.5cm,clip]{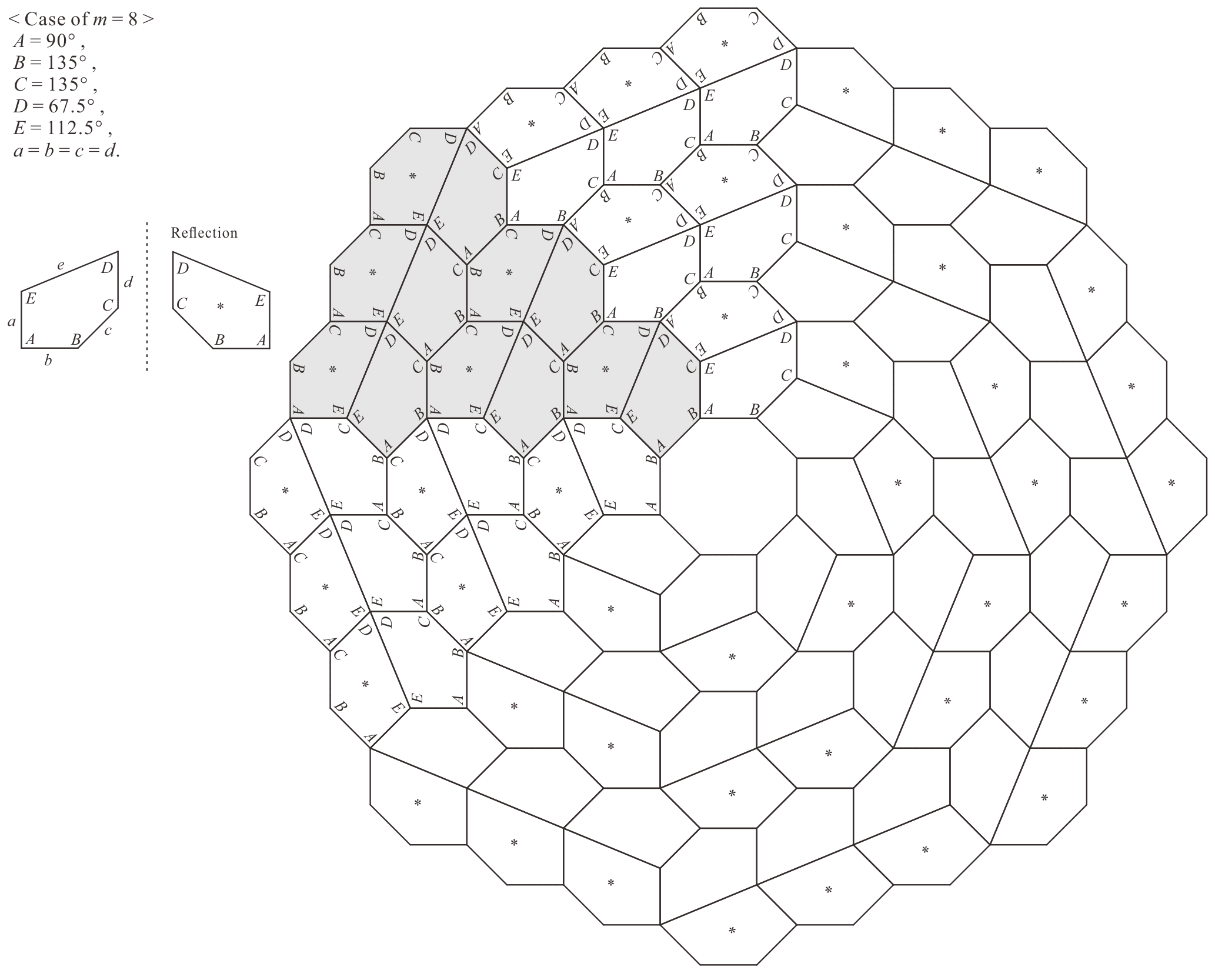} 
  \caption{{\small 
Rotationally symmetric tiling with $C_{8}$ symmetry with a regular 
convex octagonal hole at the center by a convex pentagon
} 
\label{fig17}
}
\end{figure}

\renewcommand{\figurename}{{\small Figure.}}
\begin{figure}[htbp]
 \centering\includegraphics[width=11.5cm,clip]{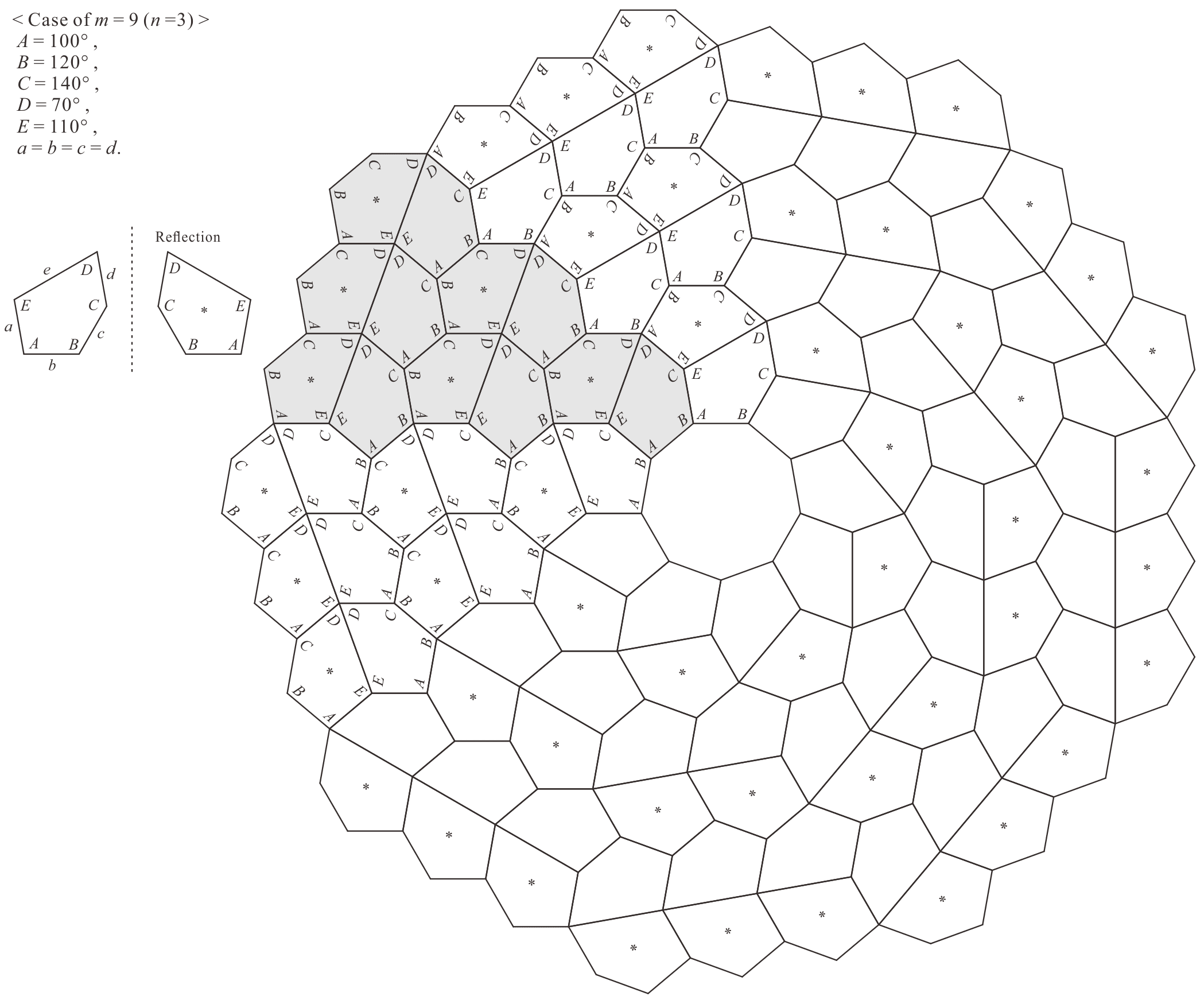} 
  \caption{{\small 
Rotationally symmetric tiling with $C_{9}$ symmetry with a regular 
convex nonagonal hole at the center by a convex pentagon
} 
\label{fig18}
}
\end{figure}

\renewcommand{\figurename}{{\small Figure.}}
\begin{figure}[htbp]
 \centering\includegraphics[width=11.5cm,clip]{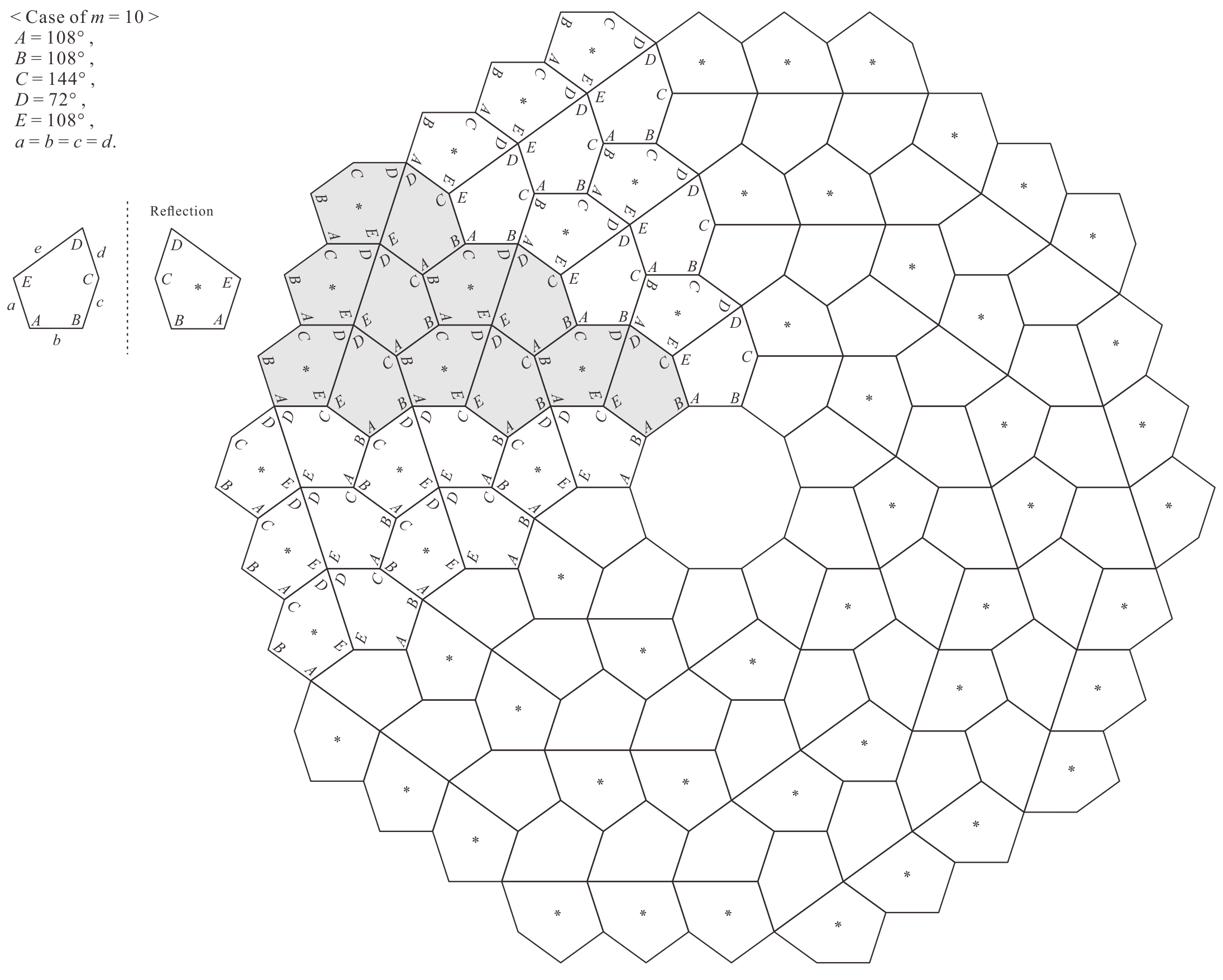} 
  \caption{{\small 
Rotationally symmetric tiling with $C_{10}$ symmetry with a regular 
convex 10-gonal hole at the center by a convex pentagon
} 
\label{fig19}
}
\end{figure}

\renewcommand{\figurename}{{\small Figure.}}
\begin{figure}[htbp]
 \centering\includegraphics[width=11.5cm,clip]{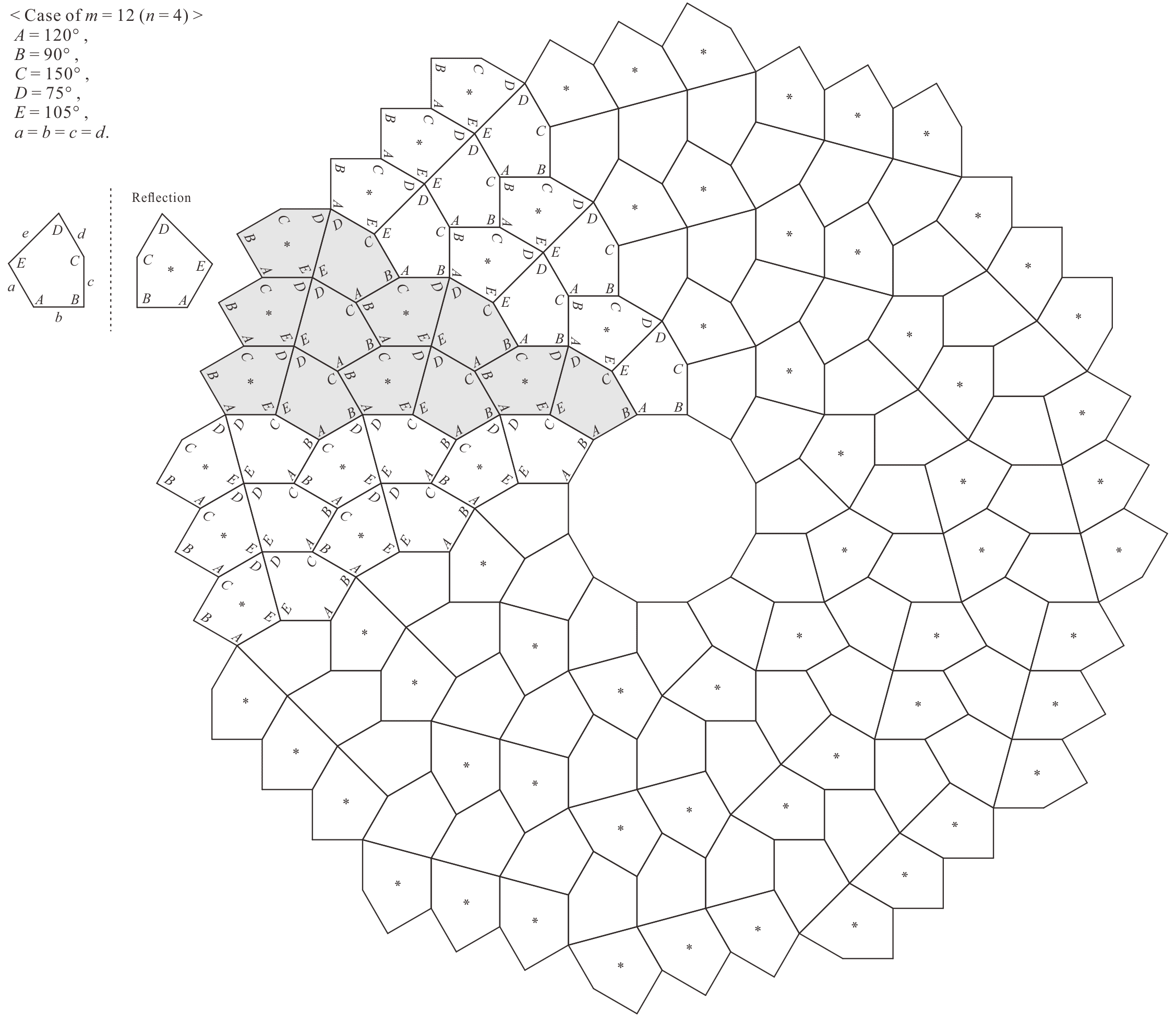} 
  \caption{{\small 
Rotationally symmetric tiling with $C_{12}$ symmetry with a regular 
convex 12-gonal hole at the center by a convex pentagon
} 
\label{fig20}
}
\end{figure}

\renewcommand{\figurename}{{\small Figure.}}
\begin{figure}[htbp]
 \centering\includegraphics[width=11.5cm,clip]{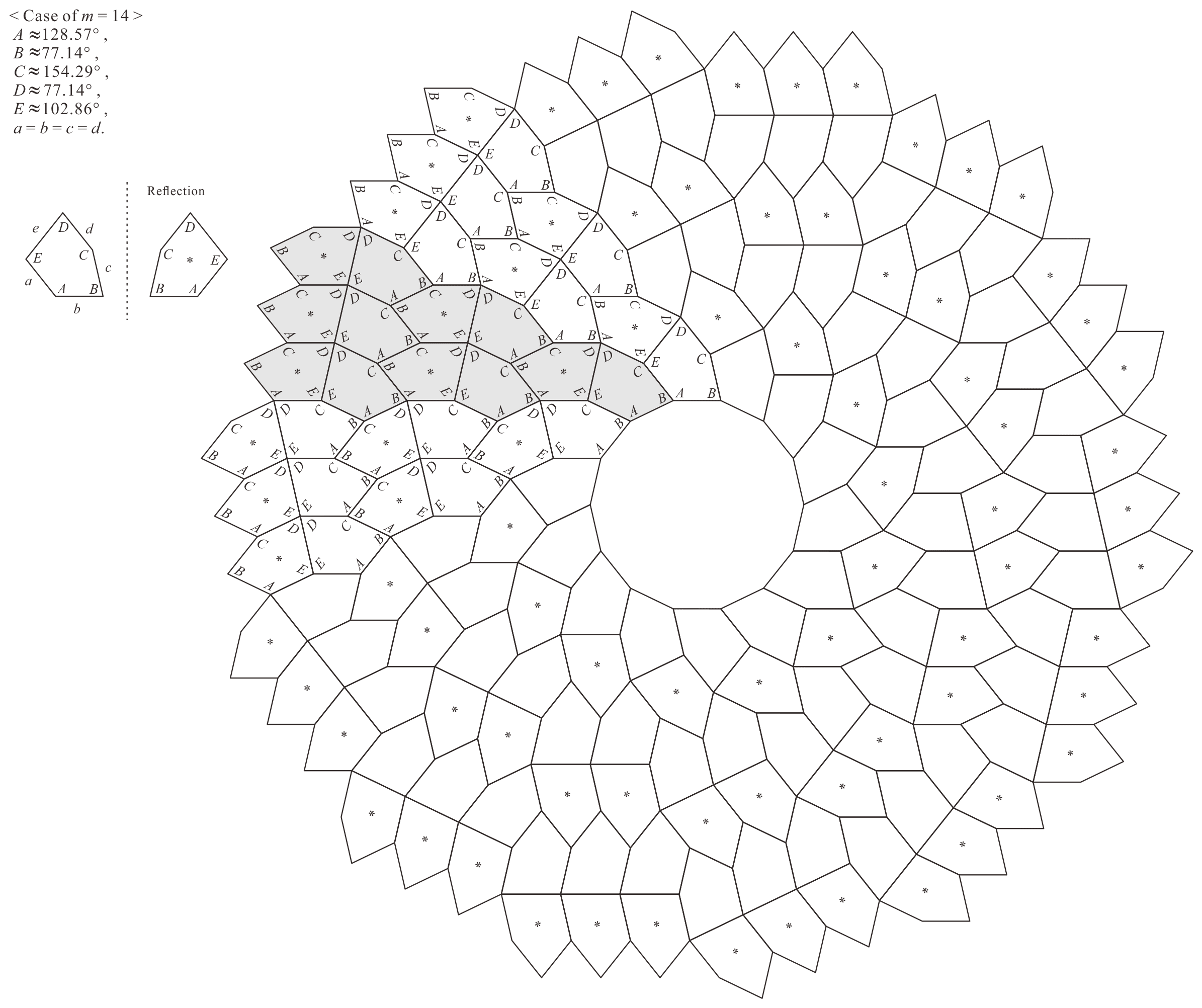} 
  \caption{{\small 
Rotationally symmetric tiling with $C_{14}$ symmetry with a regular 
convex 14-gonal hole at the center by a convex pentagon
} 
\label{fig21}
}
\end{figure}

\renewcommand{\figurename}{{\small Figure.}}
\begin{figure}[htbp]
 \centering\includegraphics[width=11.5cm,clip]{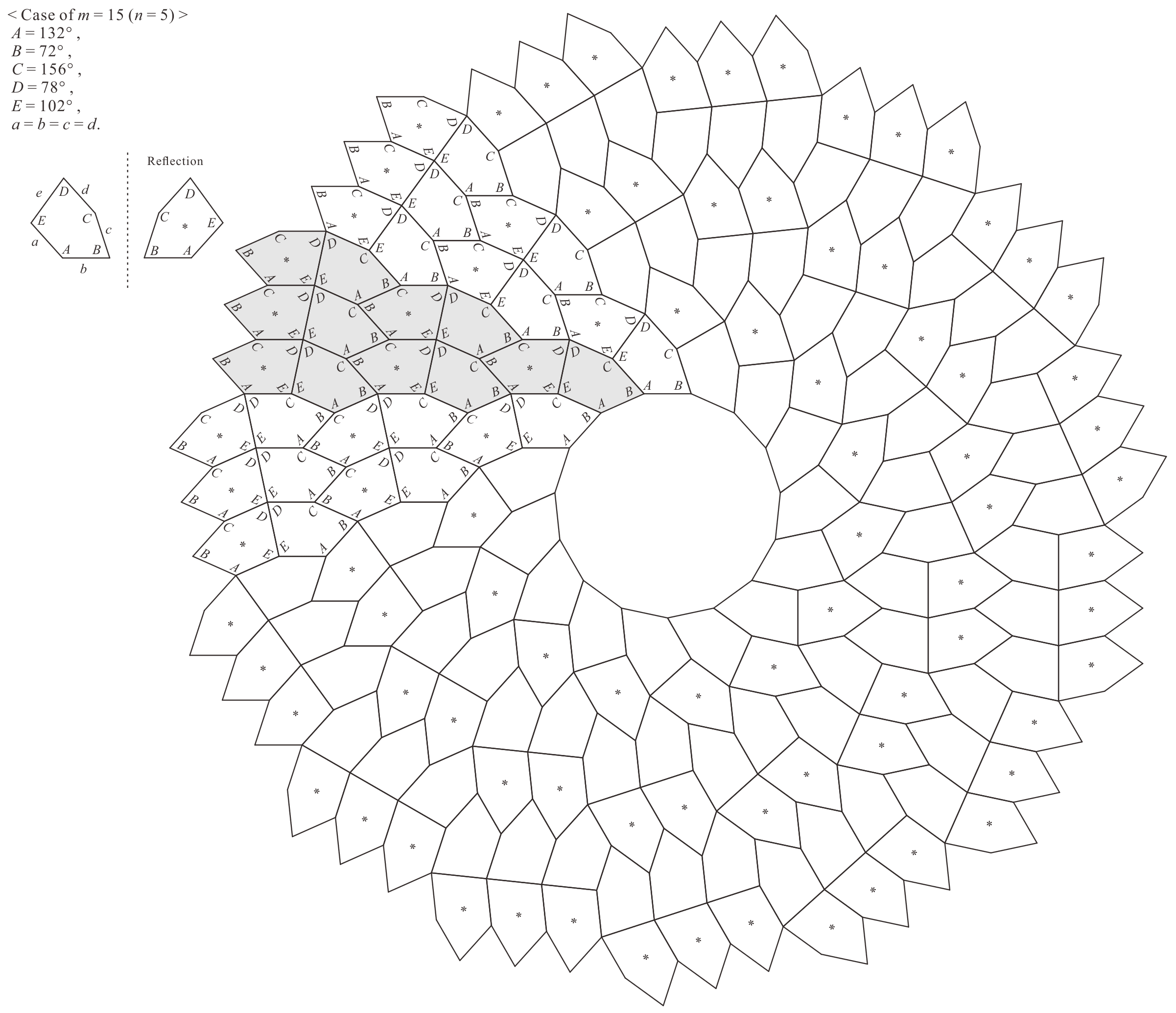} 
  \caption{{\small 
Rotationally symmetric tiling with $C_{15}$ symmetry with a regular 
convex 15-gonal hole at the center by a convex pentagon
} 
\label{fig22}
}
\end{figure}

\renewcommand{\figurename}{{\small Figure.}}
\begin{figure}[htbp]
 \centering\includegraphics[width=15cm,clip]{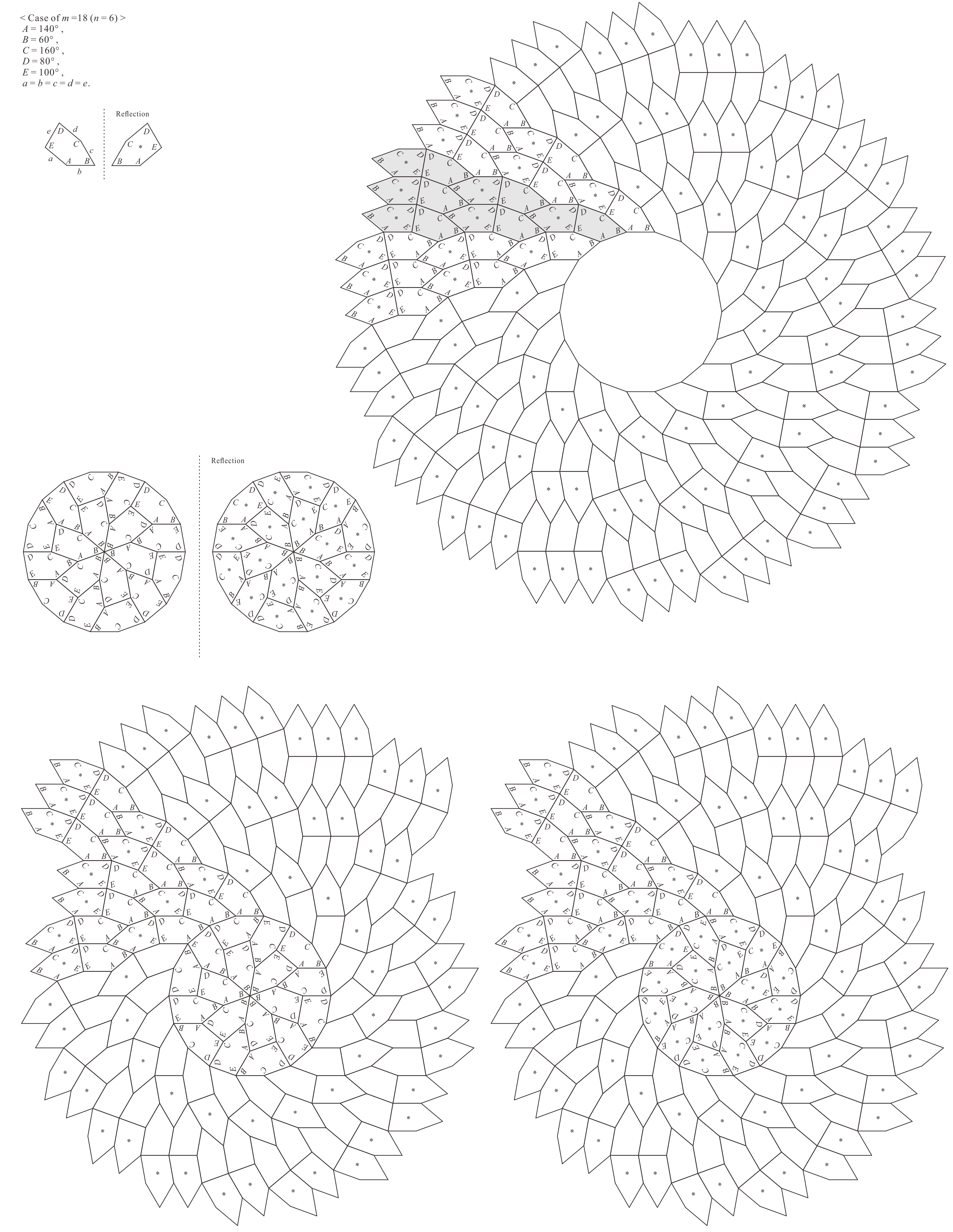} 
  \caption{{\small 
Rotationally symmetric tiling with a regular 
convex 18-gon at the center by an equilateral convex pentagon
} 
\label{fig23}
}
\end{figure}

\renewcommand{\figurename}{{\small Figure.}}
\begin{figure}[htbp]
 \centering\includegraphics[width=15cm,clip]{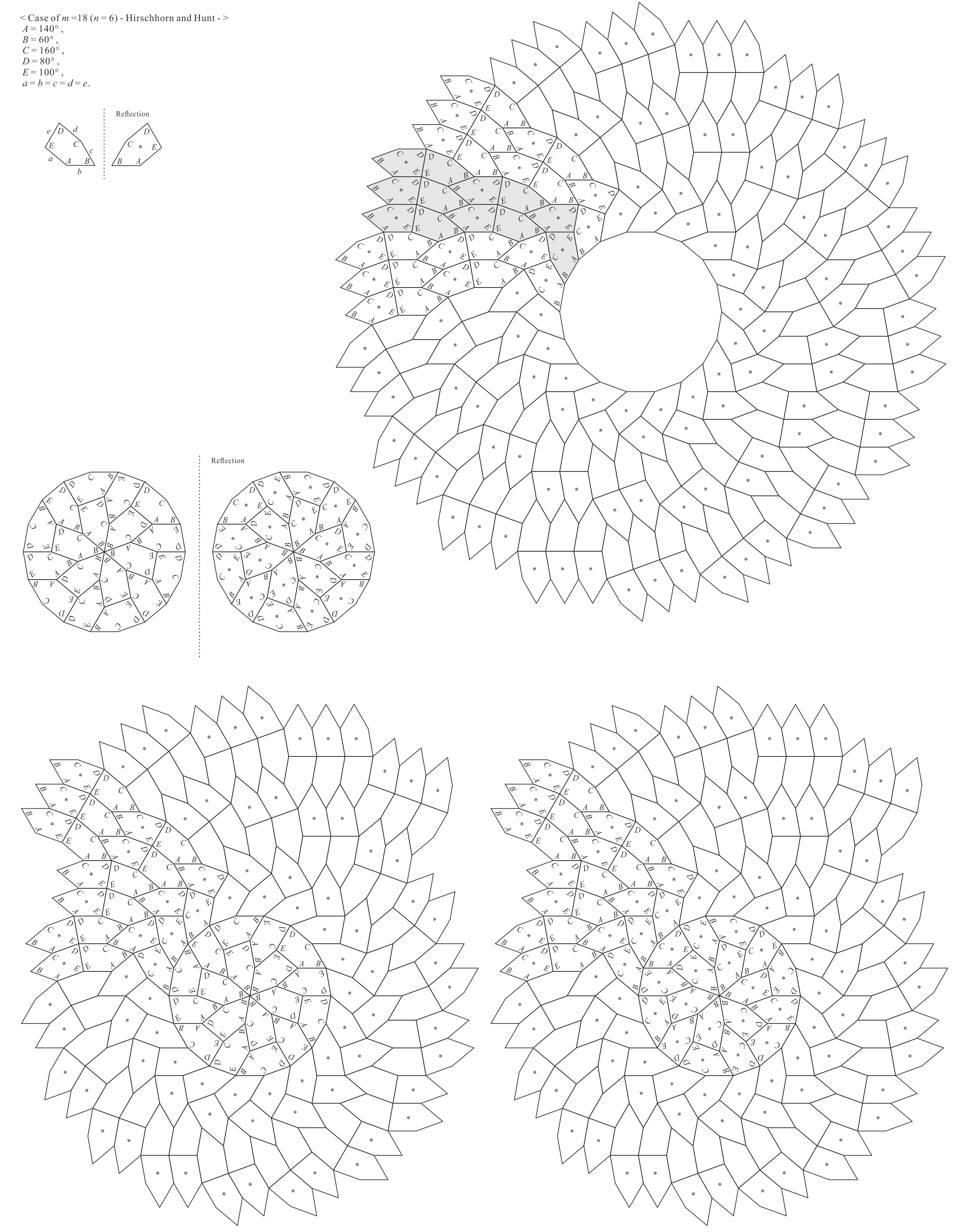} 
  \caption{{\small 
Hirschhorn and Hunt's rotationally symmetric tiling with a regular 
convex 18-gon at the center by an equilateral convex pentagon
} 
\label{fig24}
}
\end{figure}

\renewcommand{\figurename}{{\small Figure.}}
\begin{figure}[htbp]
 \centering\includegraphics[width=10.5cm,clip]{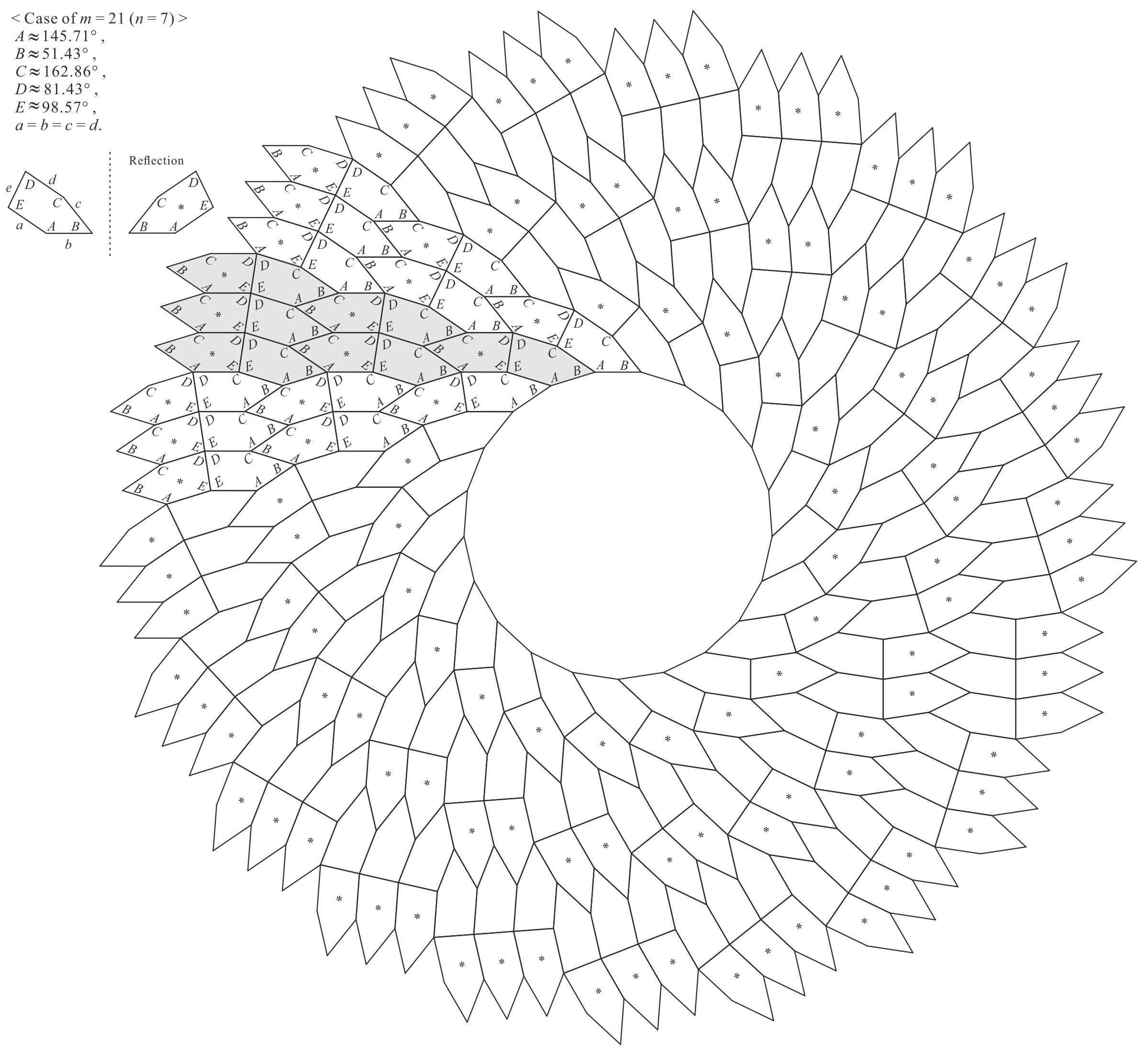} 
  \caption{{\small 
Rotationally symmetric tiling with $C_{21}$ symmetry with a regular 
convex 21-gonal hole at the center by a convex pentagon
} 
\label{fig25}
}
\end{figure}

\renewcommand{\figurename}{{\small Figure.}}
\begin{figure}[htbp]
 \centering\includegraphics[width=10.5cm,clip]{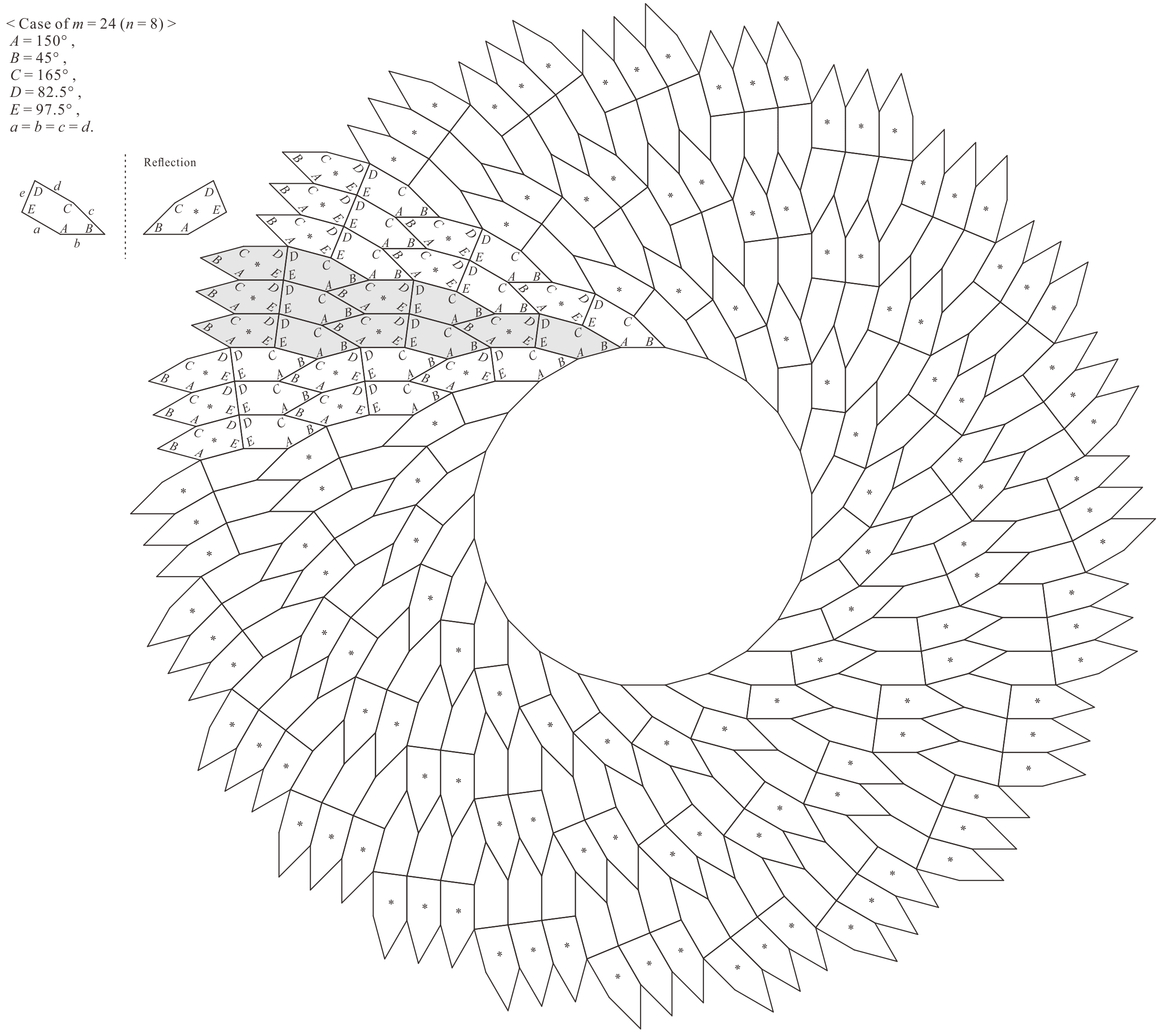} 
  \caption{{\small 
Rotationally symmetric tiling with $C_{24}$ symmetry with a regular 
convex 24-gonal hole at the center by a convex pentagon
} 
\label{fig26}
}
\end{figure}

\renewcommand{\figurename}{{\small Figure.}}
\begin{figure}[htbp]
 \centering\includegraphics[width=10.5cm,clip]{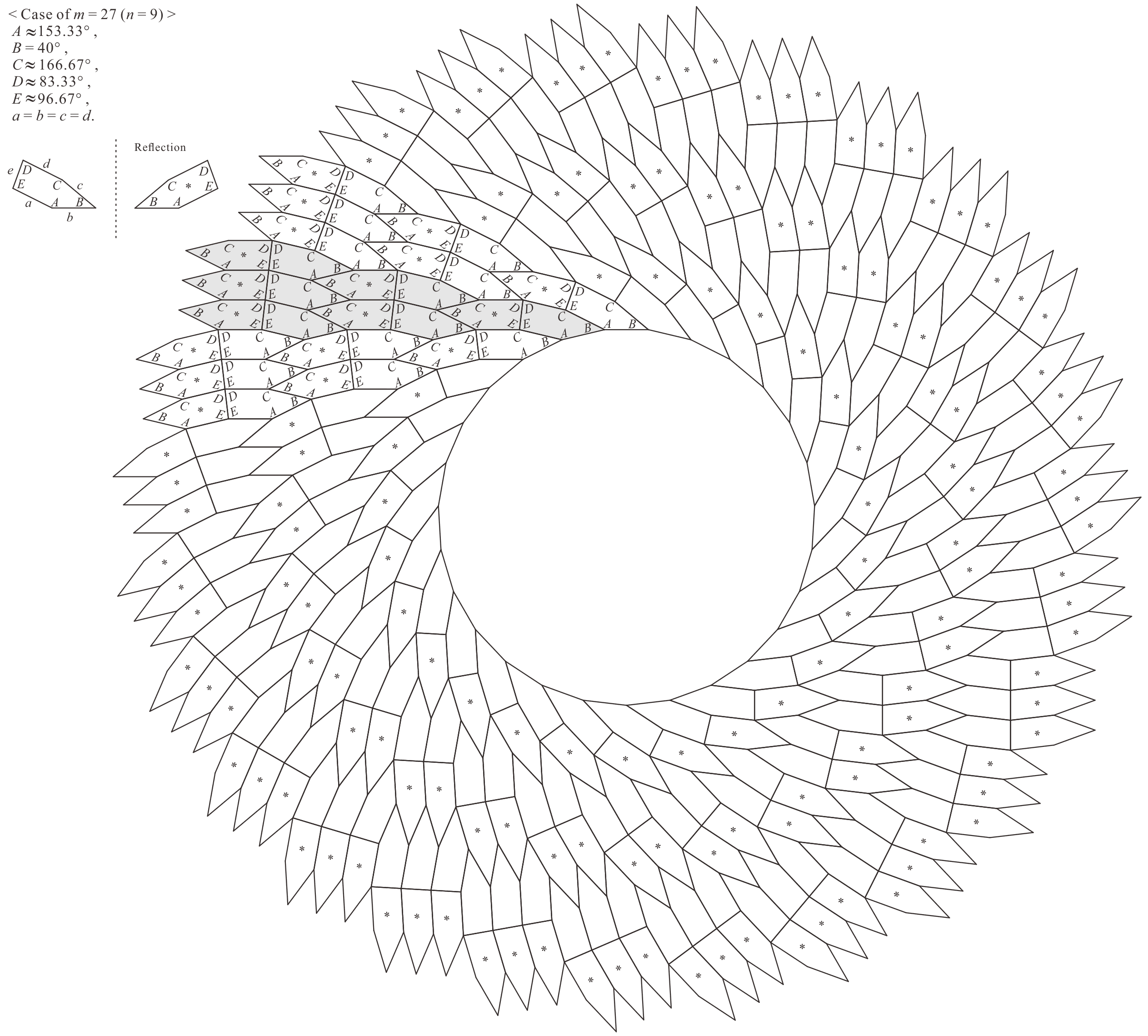} 
  \caption{{\small 
Rotationally symmetric tiling with $C_{27}$ symmetry with a regular 
convex 27-gonal hole at the center by a convex pentagon
} 
\label{fig27}
}
\end{figure}

\renewcommand{\figurename}{{\small Figure.}}
\begin{figure}[htbp]
 \centering\includegraphics[width=10cm,clip]{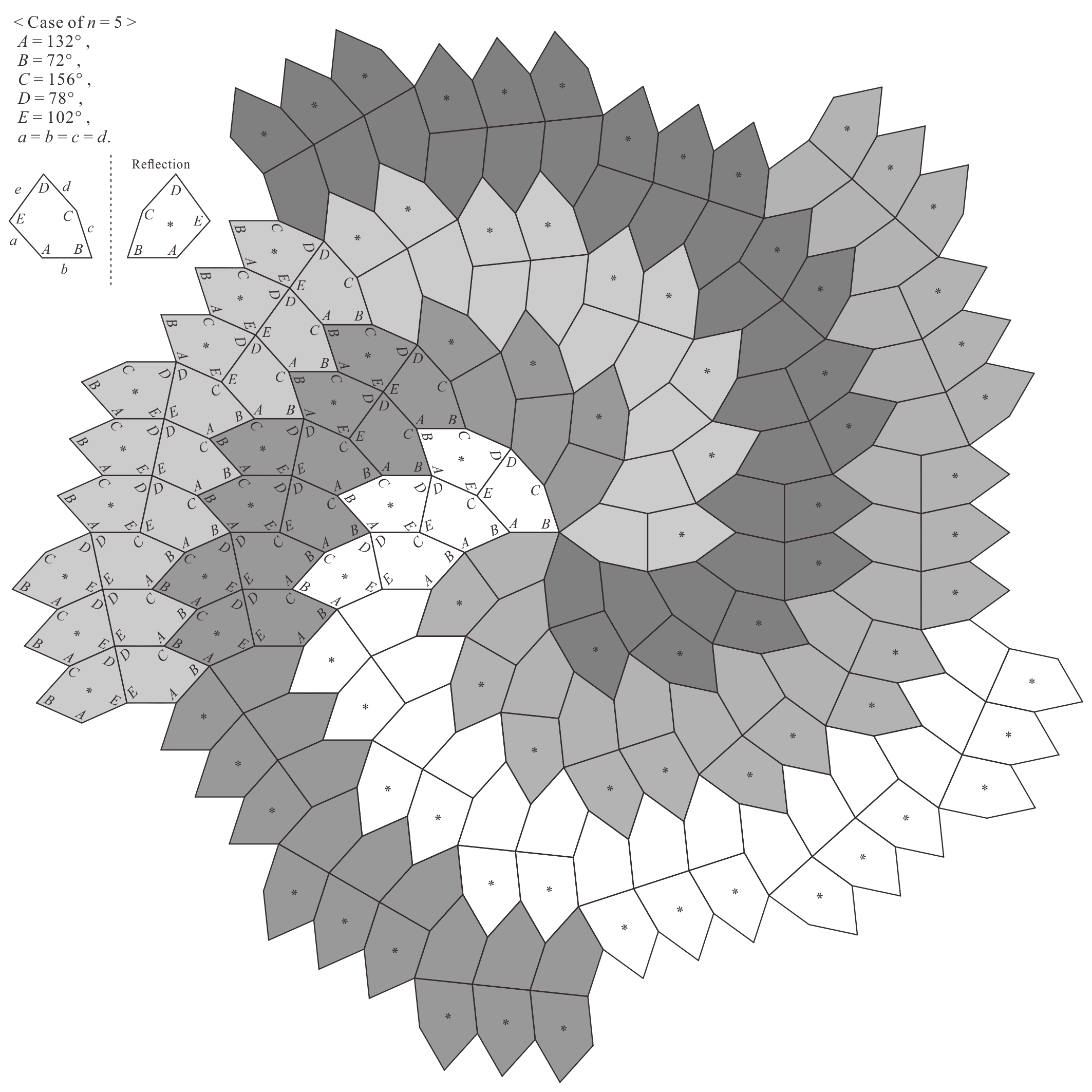} 
  \caption{{\small 
Spiral structure of five-fold rotationally symmetric edge-to-edge tiling 
by a convex pentagon
} 
\label{fig28}
}
\end{figure}

\renewcommand{\figurename}{{\small Figure.}}
\begin{figure}[htbp]
 \centering\includegraphics[width=14cm,clip]{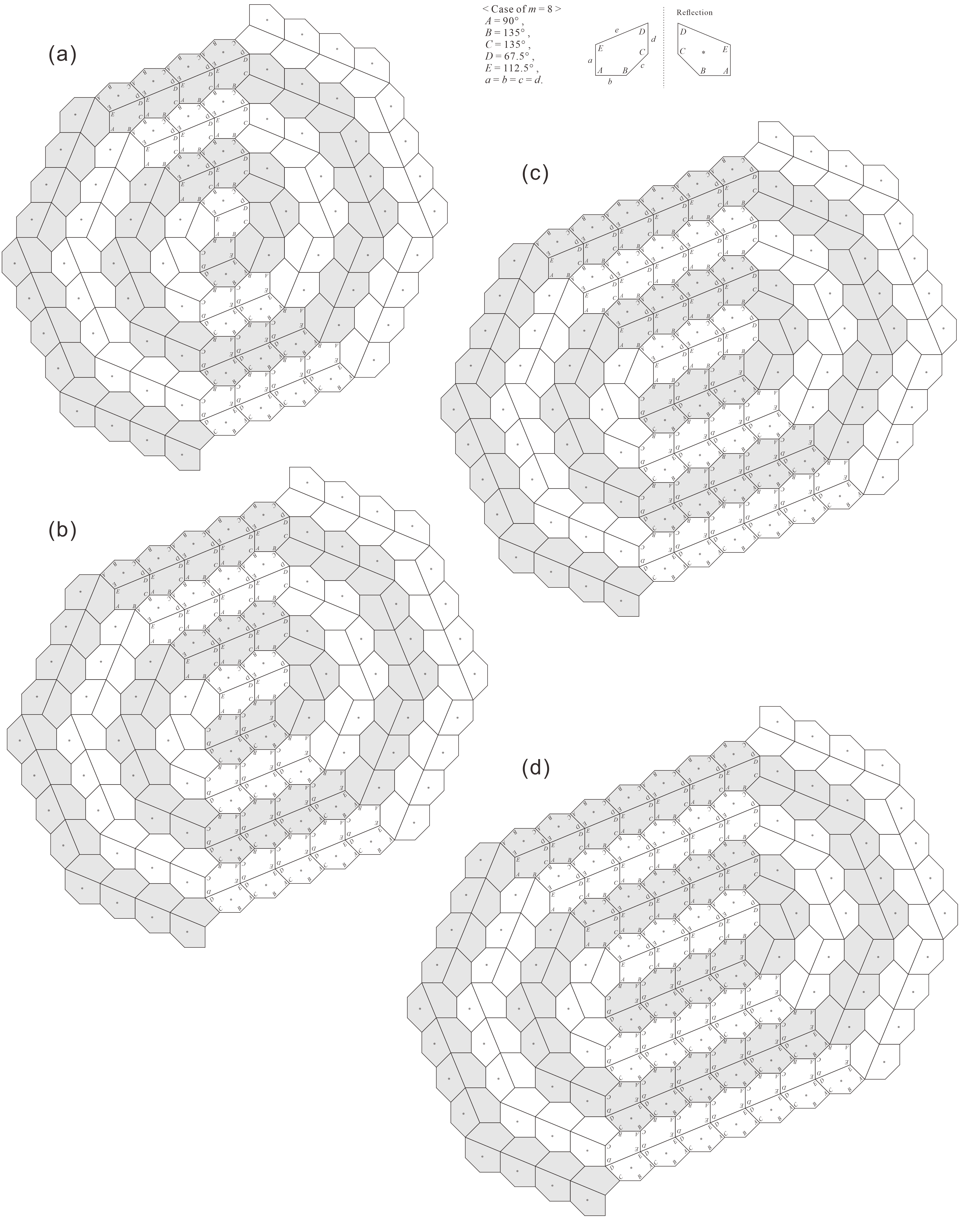} 
  \caption{{\small 
Spiral tilings with two-fold rotational symmetry by a convex pentagon 
of $m = 8$ in Table~\ref{tab2}
} 
\label{fig29}
}
\end{figure}

\renewcommand{\figurename}{{\small Figure.}}
\begin{figure}[htbp]
 \centering\includegraphics[width=15cm,clip]{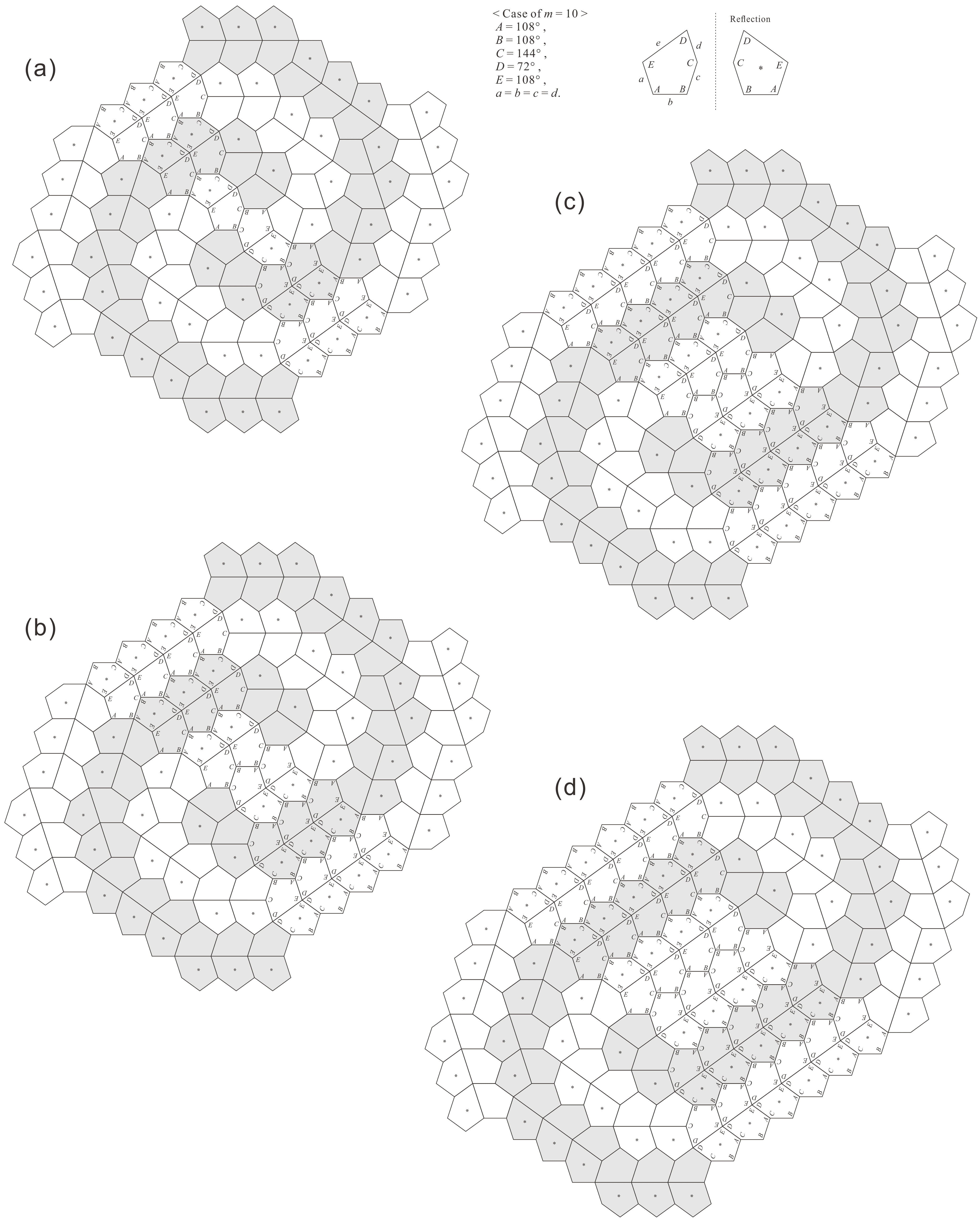} 
  \caption{{\small 
Spiral tilings with two-fold rotational symmetry by a convex pentagon 
of $m = 10$ in Table~\ref{tab2}
} 
\label{fig30}
}
\end{figure}

\renewcommand{\figurename}{{\small Figure.}}
\begin{figure}[htbp]
 \centering\includegraphics[width=15cm,clip]{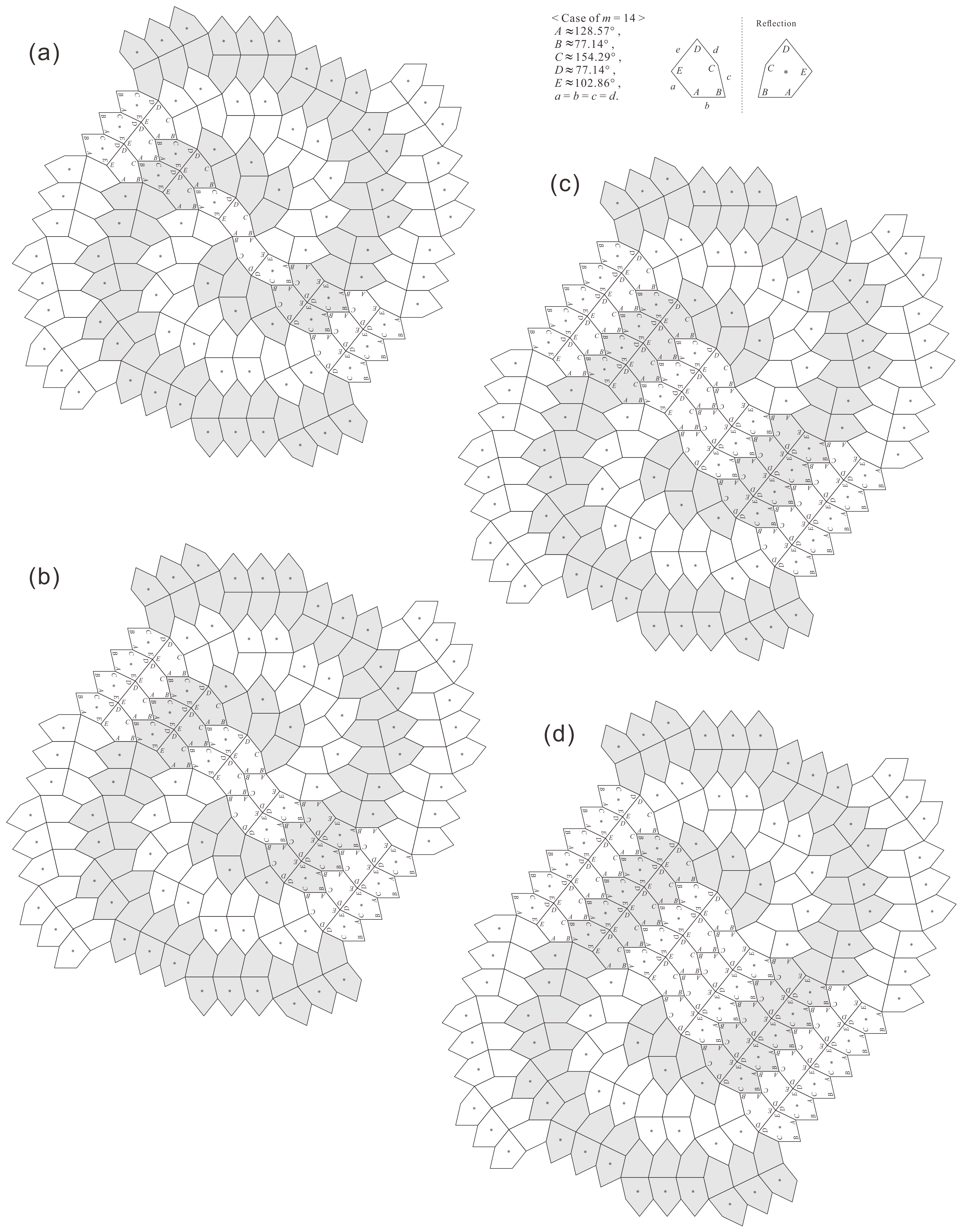} 
  \caption{{\small 
Spiral tilings with two-fold rotational symmetry by a convex pentagon 
of $m = 14$ in Table~\ref{tab2}
} 
\label{fig31}
}
\end{figure}

\renewcommand{\figurename}{{\small Figure.}}
\begin{figure}[htbp]
 \centering\includegraphics[width=15cm,clip]{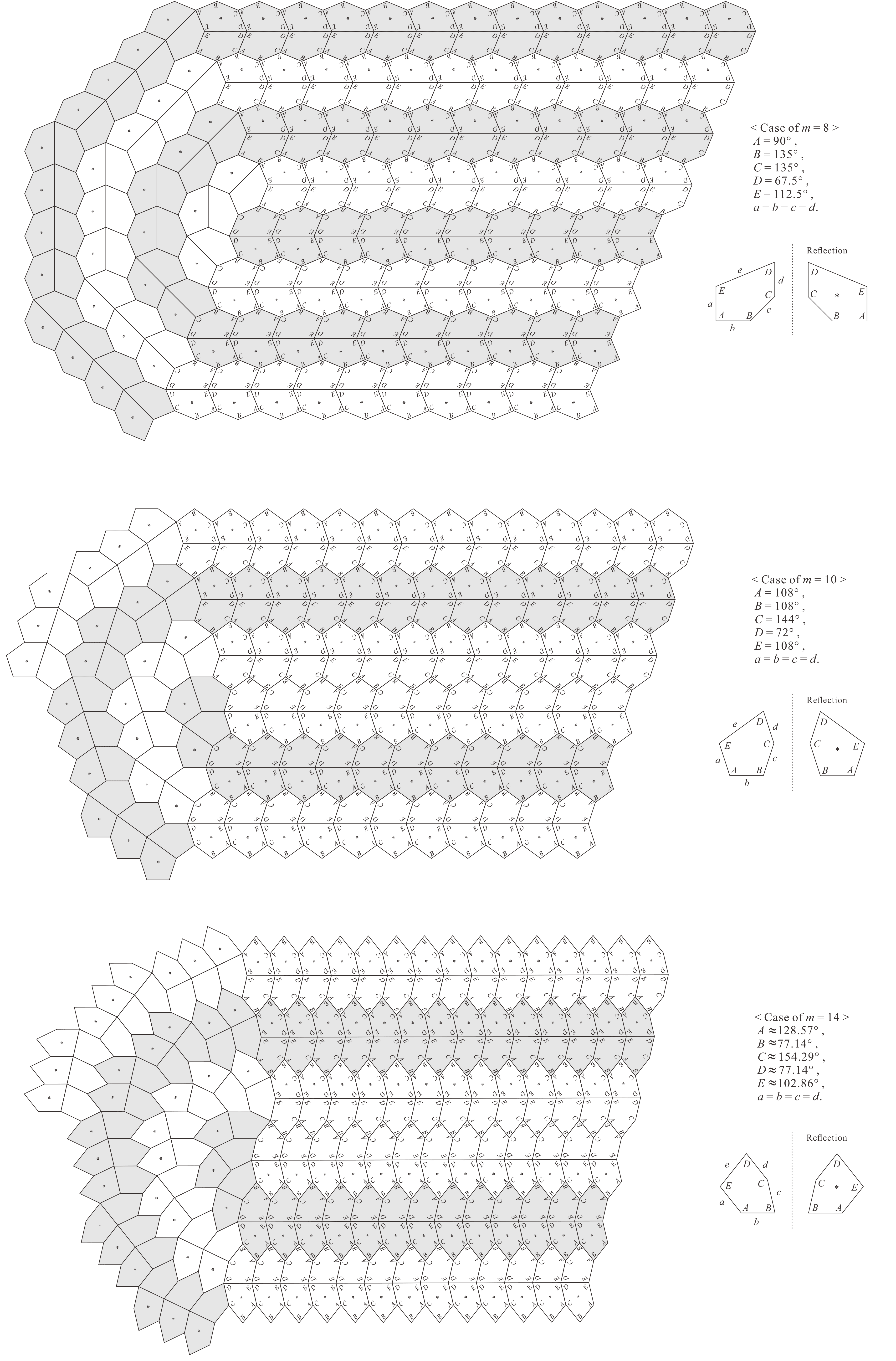} 
  \caption{{\small 
Tilings that removed one spiral structure and extend the belts of Octa-units 
using convex pentagons of $m = 8, 10, 14$ in Table~\ref{tab2}
} 
\label{fig32}
}
\end{figure}

\end{document}